\documentclass[
	english
]{scrartcl}

\usepackage[T1]{fontenc}
\usepackage[utf8]{inputenc}

\usepackage{amsmath,amsfonts,amsthm,amssymb}
\usepackage[colorlinks,citecolor=blue]{hyperref}
\usepackage{cancel}
\usepackage{relsize}
\usepackage{algorithmic}
\usepackage{booktabs,longtable}
\usepackage{rotating}

\usepackage{changes}

\usepackage{enumitem}
\setlist[enumerate]{label=(\alph*)}

\usepackage[figure]{hypcap}
\usepackage{graphicx}
\graphicspath{{./img/}}
\usepackage{subcaption}

\usepackage{preamble}

\usepackage{marginnote}

\usepackage[capitalise,nameinlink]{cleveref}
\allowdisplaybreaks

\numberwithin{equation}{section}

\crefname{assumption}{Assumption}{Assumptions}
\crefname{table}{Table}{Tables}
\Crefname{ALC@unique}{Step}{Steps}

\newcommand\norm[1]{\left\Vert#1\right\Vert}
\newcommand\nnorm[1]{\Vert#1\Vert}

\newcommand\N{\mathbb{N}}
\newcommand\R{\mathbb{R}}

\newcommand\LL{\mathcal L}

\newcommand\tto{\rightrightarrows}

\newcommand{\dom}{\operatorname{dom}}

\newcommand{\gph}{\operatorname{gph}}
\newcommand{\epi}{\operatorname{epi}}

\DeclareMathOperator*{\argmin}{\operatorname{argmin}}
\DeclareMathOperator*{\subjectto}{\operatorname{s.t.}}

\DeclareMathAlphabet{\mathpzc}{OT1}{pzc}{m}{it}

\newtheorem{theorem}{Theorem}[section]
\newtheorem{lemma}[theorem]{Lemma}
\newtheorem{proposition}[theorem]{Proposition}

\newtheorem{remark}[theorem]{Remark}

\newtheorem{example}[theorem]{Example}

\makeatletter
\long\def\@firstoffiveparen#1#2#3#4#5{\textup{\tagform@{#1}}}
\def\eqref@nolink#1{\textup{\tagform@{\ref*{#1}}}}
\def\eqref@link#1{%
\Hy@safe@activestrue
\expandafter\@setref\csname r@#1\endcsname\@firstoffiveparen{#1}%
\Hy@safe@activesfalse
}
\protected\def\eqref{\@ifstar\eqref@nolink\eqref@link}
\makeatother

\definecolor{mygreen}{rgb}{0.0,0.7,0.0}
\definecolor{mybrown}{rgb}{0.5,0.5,0.0}

\begin{document}

\title{%
	Duality-based single-level reformulations 
	of bilevel optimization problems
	}%
\author{%
	Stephan Dempe%
	\footnote{%
		Freiberg University for Mining and Technology,
		Faculty of Mathematics and Computer Science,
		09599 Freiberg,
		Germany,
		\email{dempe@extern.tu-freiberg.de},
		\orcid{0000-0001-6344-5152}
		}
	\and
	Patrick Mehlitz%
	\footnote{%
		Philipps-Universit\"at Marburg,
		Department of Mathematics and Computer Science,
		35032 Marburg,
		Germany,
		\email{mehlitz@uni-marburg.de},
		\orcid{0000-0002-9355-850X}%
		}
	}

\publishers{}
\maketitle

\begin{abstract}
	 Usually, bilevel optimization problems need to be transformed into single-level ones
	 in order to derive optimality conditions and solution algorithms.
	 Among the available approaches, the replacement of the lower-level problem
	 by means of duality relations became popular quite recently.
	 We revisit three realizations of this idea which are based
	 on the lower-level Lagrange, Wolfe, and Mond--Weir dual problem.
	 The resulting single-level surrogate problems
	 are equivalent to the original bilevel optimization problem from the viewpoint
	 of global minimizers under mild assumptions.
	 However, all these reformulations suffer from the appearance of so-called implicit
	 variables, i.e., surrogate variables which do not enter the objective function
	 but appear in the feasible set for modeling purposes.
	 Treating implicit variables as explicit ones has been shown to be problematic
	 when locally optimal solutions, stationary points, and applicable constraint
	 qualifications are compared to the original problem.
	 Indeed, we illustrate that the same difficulties have to be faced when using
	 these duality-based reformulations.
	 Furthermore, we show that the Mangasarian--Fromovitz constraint qualification is
	 likely to be violated at each feasible point of these reformulations,
	 contrasting assertions in some recently published papers.
\end{abstract}


\begin{msc}	
	\mscLink{90C30}, \mscLink{90C33}, \mscLink{90C46}
\end{msc}

\section{Introduction}\label{sec:introduction}

Bilevel optimization problems are hierarchical optimization problems
where the feasible set of the upper-level (or leader's) problem is constrained
by the graph of the solution set mapping of the parametric lower-level (or follower's) problem.
This model has been first formulated by Heinrich von Stackelberg in 1934
in the context of an economic application, see \cite{Stackelberg1934}.
About 40 years later, this problem has been introduced to the optimization community,
see \cite{BrackenMcGill1973,KornajLiptak1965}.
The recent bibliography \cite{Dempe2020} counts more than 1500 references
which are related to bilevel optimization,
including more than 65 Ph.\ D.\ theses at universities from all over the world.
Typical applications of bilevel optimization address problems in the context of
data science, energy markets, finance, and logistics.
For a detailed introduction to the topic, the interested reader is referred to the monographs 
\cite{Bard1998,Dempe2002,DempeKalashnikovPerezValdesKalashnykova2015,ShimizuIshizukaBard1997}.

Due to the implicit nature of its constraints,
the bilevel optimization problem usually needs to be transformed into a single-level one
in order to derive optimality conditions and solution algorithms.
This can be done using different approaches.
Historically, the first approach was to replace the lower-level problem using its Karush--Kuhn--Tucker (KKT) conditions.
The resulting single-level optimization problem is a so-called mathematical program with complementarity constraints (MPCC).
It is well known that the Mangasarian--Fromovitz constraint qualification (MFCQ)
is violated at all feasible points of an MPCC,
and the feasible set of this problem is of highly combinatorial structure.
Despite these difficulties, numerous contributions dealing with necessary optimality conditions and solution algorithms for MPCCs
can be found in the literature, see e.g.\ \cite{FerrisMangasarianPang2001,LuoPangRalph1996} and the references therein.
In principle, these could be used to tackle the bilevel optimization problem via the aforementioned single-level reformulation,
see \cite{KimLeyfferMunson2020} for an overview.
However, it has been shown in \cite{DempeDutta2012} that the associated MPCC,
which is referred to as the KKT reformulation of the bilevel optimization problem in the literature,
is equivalent to the original bilevel optimization problem only in terms of globally optimal solutions in general.
A locally optimal solution of the KKT reformulation does not need to be related 
to a locally optimal solution of the bilevel optimization problem.
Similar difficulties can be observed when comparing merely stationary points and the associated constraint qualifications,
see \cite{AdamHenrionOutrata2018}.
These difficulties can be traced back to the treatment of lower-level Lagrange multipliers
as explicit variables in the KKT reformulation.
Similar phenomena can be observed in many other settings of mathematical optimization
where so-called implicit variables are treated as explicit ones,
see \cite{BenkoMehlitz2021} for an overview.

A second, widely used approach goes back to Outrata, see \cite{Outrata1990}.
He suggested to represent optimality for the lower-level problem
in the upper-level constraints by adding the lower-level constraints together with an inequality
bounding the objective function value of the follower's problem 
by the associated optimal value function from above.
This results in a nonconvex, nondifferentiable optimization problem 
which is fully equivalent to the bilevel optimization problem.
Unfortunately, any nonsmooth variant of the MFCQ
is violated at all feasible points of the resulting surrogate problem,
see e.g.\ \cite[Proposition~3.2]{YeZhu1995}.
A so-called partial calmness condition has been suggested in \cite{YeZhu1995}, 
making an exact penalization approach possible
and allowing for the derivation of necessary optimality conditions,
see e.g.\ \cite{DempeDuttaMordukhovich2007,DempeZemkoho2011,MordukhovichNamPhan2012}.
It has been documented in \cite{HenrionSurowiec2011,MehlitzMinchenkoZemkoho2020}
that the partial calmness condition is comparatively restrictive.
Recently, it has been suggested to replace the lower-level value function by its Moreau envelope,
see e.g.\ \cite{BaiYeZeng2023}, but this reformulation, essentially, suffers from similar shortcomings.

A combination of the KKT and value function approaches is reasonable as well,
see \cite{YeZhu2010}, and the addition of lower-level second-order conditions
can also be taken into account, see \cite{MaYaoYeZhang2023}.
Due to the addition of superfluous constraints, partial-calmness-type conditions
are easier to hold than for the mere value function reformulation.
However, the numerical treatment of the reformulated problem is highly challenging
due to the large number of complicated constraints.
Other possible reformulation approaches use variational inequalities, see \cite{YeZhuZhu1997},
semidefinite optimization, see \cite{Bard1984},
the lower-level KKT conditions without an explicit use of a multiplier, see \cite{DempeZemkoho2013},
or the lower-level Fritz--John conditions, see \cite{AllendeStill2013}.

More recently, it has been suggested to replace the lower-level problem exploiting a suitable
(strong) duality relation. More precisely, the constraints of the lower-level problem and its
dual as well as an inequality constraint enforcing strong duality are added to the constraints
of the upper-level problem while the original lower-level problem is dropped.
This approach has been worked out based on the lower-level Wolfe dual, see \cite{Wolfe1961},
in \cite{LiLinZhangZhu2022}, and some earlier contributions already made use of this idea,
see \cite{Cerulli2021,CerulliDAmbrosioLibertiPelegrin2021,DiehlHouskaSteinSteuermann2013}.
Similarly, one can utilize the lower-level Mond--Weir dual, see \cite{MondWeir1981},
to derive a single-level reformulation of the bilevel optimization problem, see \cite{LiLinZhu2024}. 
In principle, it is also possible to exploit the lower-level Lagrange dual for this purpose,
see \cite{OuattaraAswani2018} for the analysis,
and this approach has been used to solve bilevel optimization problems with a
linear lower-level problem quite frequently.
The duality gap of the Lagrange dual associated with the lower-level problem has been penalized 
to derive a solution algorithm for linear bilevel optimization problems in 
\cite{AboussororMansouri2005,AnandalingamWhite1990,CalveteGale2004}. 
Enforcing strong duality for the linear lower-level problem, 
an evolutionary algorithm for the solution of the superordinate bilevel optimization problem
has been proposed in \cite{LiFang2012}.
In \cite{TuyGhannadan1998}, the authors exploit the lower-level Lagrange dual to
transform a linear bilevel optimization problem into a reverse convex program 
which can be solved globally by respective algorithms.
Application-driven publications which exploit this approach are, exemplary,
\cite{GarciaHerrerosMisraArslanMehtaGrossmann2015,TangSunHauser2020}.

In this survey paper, we review the transformation approaches which are based on the lower-level
Lagrange, Wolfe, and Mond--Weir dual problem and provide a thorough analysis of the
relationship between the original bilevel optimization problem and its respective
single-level reformulation. In this regard, we complement the results which have been
stated in \cite{LiLinZhangZhu2022,LiLinZhu2024,OuattaraAswani2018}.
We point out equivalence between the bilevel optimization problem and its transformation 
when globally optimal solutions are investigated. 
However, we also show that locally optimal solutions of all three single-level reformulations 
are, in general, not related to locally optimal solutions of the original bilevel optimization problem,
and this issue is illustrated by means of an example in all three cases.
As all these single-level reformulations suffer from the presence of artificial variables
which are mainly introduced for modeling purposes, this phenomenon is not surprising in the light
of \cite{BenkoMehlitz2021} where issues related to implicit variables, which are treated as explicit
ones, are studied in rigorous detail.
In \cite{LiLinZhangZhu2022,LiLinZhu2024}, the authors claim that MFCQ
can hold at certain feasible points of the single-level reformulations which are based on the lower-level
Wolfe and Mond--Weir dual. Here, we refute these assertions.
Indeed, we prove that this qualification condition fails to hold at each feasible point of those single-level reformulations.
Furthermore, we show that a suitable nonsmooth version of MFCQ is likely to fail
at the feasible points of the single-level reformulation which exploits the lower-level Lagrange dual problem
under mild assumptions. This way, it is illustrated that lower-level duality-based reformulations of the bilevel
optimization problem suffer from the same intrinsic difficulties as the classical KKT reformulation.

The remainder of the paper is organized as follows.
After a brief discussion of the exploited notation in \cref{sec:notation},
three duality approaches for nonlinear optimization problems (Lagrange, Wolfe, and Mond--Weir duality) 
are reviewed in \cref{sec:duality}. 
In \cref{sec:bpp_intro}, we formally introduce the bilevel optimization problem and recall its standard transformations.
\cref{sec:duality_based_trafos} comprises the main results of this paper.
We study the single-level reformulations of the bilevel optimization problem which are based on the lower-level
Lagrange, Wolfe, and Mond--Weir dual in \cref{sec:trafo_Lagrange_duality,sec:trafo_Wolfe_duality,sec:trafo_Mond_Weir_duality},
respectively. 
A brief comparison is provided in \cref{sec:comparison}. 
Some concluding remarks close the paper in \cref{sec:conclusions}.

\section{Preliminaries}

In this section, we first comment on the notation, which is exploited in this manuscript, in \cref{sec:notation}.
Afterwards, the essentials of Lagrange, Wolfe, and Mond--Weir duality
for convex optimization problems are briefly recalled in \cref{sec:duality}.

\subsection{Notation}\label{sec:notation}

Let $\N$ and $\N_0$ be the sets of all positive and nonnegative integers, respectively.
Furthermore, $\R$, $\R_+$, and $\R_-$ are the sets of all real numbers, nonnegative real numbers,
and nonpositive real numbers, respectively,
and $\overline\R:=\R\cup\{-\infty,\infty\}$ is the extended real line.
The Euclidean norm of a vector $w\in\R^s$ will be represented by $\norm{w}$.

For differentiable functions $q\colon\R^s\to\R^t$ as well as $p\colon\R^s\to\R$
and some point $w\in\R^s$,
$q'(w)\in\R^{t\times s}$ denotes the Jacobian of $q$ at $w$
and $\nabla p(w):=p'(w)^\top\in\R^s$ represents the gradient of $p$ at $w$, respectively.
Whenever $p$ is even twice differentiable at $w$,
$\nabla^2p(w):=(\nabla p)'(w)\in\R^{s\times s}$ is the Hessian of $p$ at $w$.
Partial derivatives with respect to (w.r.t.) certain blocks of variables are indicated
with the aid of indices in canonical way.

Recall that, for a closed set $\Omega\subset\R^s$ and some point $w\in\Omega$,
\[
	N_\Omega(w)
	:=
	\left\{
		\eta\in\R^s\,\middle|\,
		\begin{aligned}
			&\exists\{w^k\}_{k\in\N},\{v^k\}_{k\in\N}\subset\R^s,\,
				\exists\{\alpha_k\}_{k\in\N}\subset(0,\infty)\colon
				\\
			&\quad w^k\to w,\,\alpha_k(w^k-v^k)\to\eta,\,v^k\in\Pi_\Omega(w^k)\,\forall k\in\N
		\end{aligned}
	\right\}
\]
is called the limiting (or Mordukhovich) normal cone to $\Omega$ at $w$.
Above, for each $w'\in\R^s$,
\[
	\Pi_\Omega(w')
	:=
	\argmin\limits_v\{\nnorm{v-w'}\,|\,v\in\Omega\}
\]
is the projector of $w'$ onto $\Omega$.
For each $w'\notin\Omega$, we set $N_\Omega(w'):=\emptyset$.
Recall that, for some lower semicontinuous function $p\colon\R^s\to\overline{\R}$,
its epigraph
\[
	\epi p
	:=
	\{(w,\alpha)\in\R^s\times\R\,|\,\alpha\geq p(w)\}
\]
is closed. Hence, fixing $w\in\dom p:=\{w\in\R^s\,|\,|p(w)|<\infty\}$,
the constructions
\begin{align*}
	\partial p(w)
	&:=
	\{\eta\in\R^s\,|\,(\eta,-1)\in N_{\epi p}(w,p(w))\},
	\\
	\partial^\infty p(w)
	&:=
	\{\eta\in\R^s\,|\,(\eta,0)\in N_{\epi p}(w,p(w))\},
\end{align*}
which are called the limiting and singular subdifferential of $p$ at $w$, respectively,
are reasonable.
Let us mention that $\partial^\infty p(w)=\{0\}$ holds if and only if $p$ is
locally Lipschitz continuous at $w$, see e.g.\ \cite[Theorem~1.22]{Mordukhovich2018}.
For a closed set $\Omega\subset\R^s$, $\delta_\Omega\colon\R^s\to\overline{\R}$
denotes the indicator function of $\Omega$.
The latter vanishes on $\Omega$ and is set to $\infty$ on $\R^s\setminus\Omega$.
Hence, $\epi\delta_\Omega=\Omega\times[0,\infty)$,
and, for each $w\in\Omega$,
\[
	\partial \delta_\Omega(w)
	=
	\partial^\infty\delta_\Omega(w)
	=
	N_\Omega(w)
\]
by the product rule for limiting normals,
see e.g.\ \cite[Proposition~1.4]{Mordukhovich2018}.

For lower semicontinuous functions $q_1,\ldots,q_\nu\colon\R^s\to\overline{\R}$,
$\nu\in\N_0$,
and continuously differentiable functions $q_{\nu+1},\ldots,q_t\colon\R^s\to\R$,
$t\in\N$, $t\geq\nu$,
we are interested in the constraint system
\begin{equation}\label{eq:constraint_system}\tag{CS}
	q_i(w)\leq 0\quad\forall i\in\{1,\ldots,\nu\},
	\qquad
	q_i(w)=0\quad \forall i\in\{\nu+1,\ldots,t\}.
\end{equation}
Pick a point $\bar w\in\R^s$ which is feasible to \eqref{eq:constraint_system}
and satisfies $\bar w\in\bigcap_{i=1}^\nu\dom q_i$.
We say that the nonsmooth Mangasarian--Fromovitz constraint qualification
(NSMFCQ for brevity)
is valid at $\bar w$ whenever the system
\begin{align*}
	-\sum_{i=\nu+1}^t v_i\nabla q_i(\bar w)
	\in
	\sum_{i=1}^\nu v_i\partial q_i(\bar w)&,
	\\
	v_i\geq 0,\,v_i\,q_i(\bar w)=0
	&
	\quad
	\forall i\in\{1,\ldots,\nu\}
\end{align*}
possesses only the trivial solution $v_1=\ldots=v_t=0$.
Clearly, whenever $q_1,\ldots,q_\nu$ are continuously differentiable at $\bar w$,
	which we will assume in the remainder of this subsection,
then NSMFCQ reduces to the well-known Mangasarian--Fromovitz constraint qualification
from standard nonlinear optimization,
and for clarity of presentation, we will use the abbreviation MFCQ in this situation.
Let us recall the well-known fact that the presence of MFCQ 
implies validity of Guignard's constraint qualification (GCQ for brevity)
which demands that
\[
	T_{\mathcal Q}(\bar w)
	=
	\left\{
		d\in\R^s\,\middle|\,
		\begin{aligned}
			q'_i(\bar w)d&\leq 0&&\forall i\in\{j\in\{1,\ldots,\nu\}\,|\,q_j(\bar w)=0\}
			\\
			q'_i(\bar w)d&=0&&\forall i\in\{\nu+1,\ldots,t\}
		\end{aligned}
	\right\},
\]
where $\mathcal Q\subset\R^s$ denotes the set of all points
which satisfy the system \eqref{eq:constraint_system},
and $T_{\mathcal Q}(\bar w)$ is the standard tangent cone to $\mathcal Q$ at $\bar w$
which is defined as stated below:
\[
	T_{\mathcal Q}(\bar w)
	:=
	\left\{
		d\in\R^s\,\middle|\,
		\begin{aligned}
			&\exists\{w^k\}_{k\in\N}\subset\mathcal Q,\,
				\exists\{t_k\}_{k\in\N}\subset(0,\infty)\colon
				\\
			&\quad w^k\to \bar w,\,t_k\downarrow 0,\,(w^k-\bar w)/t_k\to d
		\end{aligned}
	\right\}.
\]
It is well known that for each local minimizer $\bar w\in\R^s$
of $\min_w\{p(w)\,|\,w\in\mathcal Q\}$,
where $p\colon\R^s\to\R$ is a continuously differentiable function,
such that GCQ holds at $\bar w$,
there exist so-called Lagrange multipliers $v\in\R^t$ such that
\[
	\begin{aligned}
		\nabla p(\bar w) + \sum_{i=1}^tv_i\nabla q_i(\bar w)=0&,
		\\
		q_i(\bar w)\leq 0,\,v_i\geq 0,\,v_i\,q_i(\bar w)=0&\quad
		\forall i\in\{1,\ldots,\nu\},
		\\
		q_i(\bar w)=0&\quad
		\forall i\in\{\nu+1,\ldots,t\}.
	\end{aligned}
\]
These conditions are referred to as the Karush--Kuhn--Tucker (KKT for brevity) conditions,
and a tuple $(\bar w,v)$ fulfilling them is called a KKT point while $\bar w$ is referred
to as stationary in this situation.
We note that, whenever $p,q_1,\ldots,q_\nu$ are convex and $q_{\nu+1},\ldots,q_t$ are affine,
then a stationary point is also a minimizer
whether a constraint qualification is valid or not,
see \cite[Theorem~3.27]{Ruszczynski2011}.

Finally, for a set-valued mapping $\Upsilon\colon\R^s\tto\R^t$,
we exploit
$\dom\Upsilon:=\{w\in\R^s\,|\,\Upsilon(w)\neq\emptyset\}$
and
$\gph\Upsilon:=\{(w,v)\in\R^s\times\R^t\,|\,v\in\Upsilon(w)\}$
in order to denote the domain and the graph of $\Upsilon$.
Let us fix $\bar w\in\dom\Upsilon$.
Then $\Upsilon$ is said to be inner semicompact at $\bar w$ w.r.t.\ its domain
whenever for each sequence $\{w^k\}_{k\in\N}\subset\dom\Upsilon$ such that $w^k\to\bar w$,
there is a sequence $\{v^k\}_{k\in\N}$ which, on the one hand,
satisfies $v^k\in\Upsilon(w^k)$ for all $k\in\N$ and, on the other hand,
possesses a bounded subsequence.
	It is immediately clear from the definition that
	whenever $\Upsilon$ possesses uniformly bounded image sets
	in a neighborhood of $\bar w$,
	then it is inner semicompact w.r.t.\ its domain at $\bar w$.
	This well-known fact will be used frequently in this paper.

\subsection{Duality in nonlinear convex optimization}\label{sec:duality}

Throughout the section, we consider the nonlinear optimization problem
\begin{equation}\label{eq:NLP}\tag{P}
	\min\limits_w\{p(w)\,|\,q(w)\leq 0\}
\end{equation}
where $p\colon\R^s\to\R$ and $q\colon\R^s\to\R^t$ are twice continuously differentiable functions
such that $p$ and the component functions $q_1,\ldots,q_t$ of $q$ are convex.
	Let us emphasize that convexity of the data functions in \eqref{eq:NLP}
	guarantees that local and global minimizers of this problem coincide.
	We will, thus, refer to the global minimizers of the convex problem \eqref{eq:NLP}
	just as minimizers.
	Similarly, when considering the maximization of a concave function
	over a convex feasible set, we will simply refer to global maximizers
	as maximizers.
The Lagrangian function $\LL\colon\R^s\times\R^t\to\R$ of \eqref{eq:NLP} is given by
\[
	\forall w\in\R^s,\,\forall v\in\R^t\colon\quad
	\LL(w,v):=p(w)+v^\top q(w).
\]
Note that $v$ plays the role of the Lagrange multiplier in this subsection.
Furthermore, we refer to
\[
	V(w):=\{v\in\R^t\,|\,\nabla_1\LL(w,v)=0,\,0\leq v\perp -q(w)\geq 0\}
\]
as the set of Lagrange multipliers associated with $w$ for \eqref{eq:NLP}.
Recall that a pair $(\bar w,\bar v)\in\R^s\times\R^t_+$ is said to be a saddle point of $\LL$
if and only if
\[
	\forall w\in\R^s,\,\forall v\in\R^t_+\colon\quad
	\LL(\bar w,v)\leq\LL(\bar w,\bar v)\leq\LL(w,\bar v),
\]
and the latter is equivalent to $\bar v\in V(\bar w)$ due to \cite[Theorem~4.7]{Ruszczynski2011}.

Recall that \eqref{eq:NLP} satisfies Slater's qualification condition whenever there exists
a point $\tilde w\in\R^s$ such that $q_i(\tilde w)<0$ holds for all $i\in\{1,\ldots,t\}$.
In this case, MFCQ and, thus, GCQ hold at each feasible point of \eqref{eq:NLP}.

\subsubsection{Lagrange duality}\label{sec:Lagrange_duality}

To start, let us define the Lagrange value function
$\varphi_\ell\colon\R^s\to\overline{\R}$ of \eqref{eq:NLP}
by means of
\[
	\forall v\in\R^t\colon\quad
	\varphi_\ell(v)
	:=
	\begin{cases}
		\inf_w\{\LL(w,v)\,|\,w\in\R^s\}	&	v\in\R^t_+,
		\\
		-\infty							&	\text{otherwise.}
	\end{cases}
\]
We note that $\varphi_\ell$ is a concave function which does not take value $\infty$.
Furthermore, since $(w,v)\mapsto \LL(w,v)-\delta_{\R^t_+}(v)$ is upper semicontinuous
by closedness of $\R^t_+$ while $\varphi_\ell$ admits the representation
\[
	\forall v\in\R^t\colon\quad
	\varphi_\ell(v)=\inf\limits_w\{\LL(w,v)-\delta_{\R^t_+}(v)\,|\,w\in\R^s\},
\]
$\varphi_\ell$ is an upper semicontinuous function as well,
see e.g.\ \cite[Theorem~4.2.2.(1)]{BankGuddatKlatteKummerTammer1982}.
The Lagrange dual problem of \eqref{eq:NLP} is then given by means of
\begin{equation}\label{eq:D_NLP_Lagrange}\tag{D$_\ell$}
	\max\limits_{\hat v}\{\varphi_\ell(\hat v)\,|\,\hat v\geq 0\}.
\end{equation}
For feasible points $w\in\R^s$ of \eqref{eq:NLP} and $\hat v\in\R^t$ of \eqref{eq:D_NLP_Lagrange},
we find
\[
	p(w)
	\geq
	\LL(w,\hat v)
	\geq
	\varphi_\ell(\hat v)
\]
by $q(w)\leq 0$ as well as $\hat v\geq 0$ and the definition of $\varphi_\ell$,
i.e., a weak duality relation holds between \eqref{eq:NLP} and \eqref{eq:D_NLP_Lagrange}.
If a constraint qualification is valid, strong duality can be established.
In the subsequently stated proposition, we summarize the essential duality relations
between \eqref{eq:NLP} and \eqref{eq:D_NLP_Lagrange}.

\begin{proposition}\label{prop:Lagrange_duality}
	\begin{enumerate}
	\item\label{item:Lagrange_weak}
		For feasible points $w\in\R^s$ of \eqref{eq:NLP}
		and $\hat v\in\R^t$ of \eqref{eq:D_NLP_Lagrange},
		we have $p(w)\geq\varphi_\ell(\hat v)$.
		Particularly, the infimal value of \eqref{eq:NLP} is not smaller
		than the supremal value of \eqref{eq:D_NLP_Lagrange}.
	\item\label{item:Lagrange_strong}
		Let $w\in\R^s$ be a minimizer of \eqref{eq:NLP}
		where GCQ holds.
		Then there exists $v\in\R^t$ such that $v$
		is a maximizer of \eqref{eq:D_NLP_Lagrange}
		and $p(w)=\varphi_\ell(v)$.
	\item\label{item:Lagrange_saddle_point}
		Fix $(\bar w,\bar v)\in\R^s\times\R^t_+$.
		Then the following statements are equivalent.
		\begin{itemize}
			\item The pair $(\bar w,\bar v)$ is a saddle point of $\LL$.
			\item We have $\bar v\in V(\bar w)$.
			\item The point $\bar w$ is feasible to \eqref{eq:NLP},
				the point $\bar v$ is feasible to \eqref{eq:D_NLP_Lagrange},
				and $p(\bar w)=\varphi_\ell(\bar v)$.
		\end{itemize}
	\end{enumerate}
\end{proposition}
\begin{proof}
	Assertion~\ref{item:Lagrange_weak} has been validated above.
	For the proof of assertion~\ref{item:Lagrange_strong},
	observe that the assumptions guarantee the existence of $v\in V(w)$.
	Thus, $(w,v)$ is a saddle
	point of $\LL$ which particularly yields
	\[
		\forall\hat w\in\R^s\colon\quad
		\LL(w,v)\leq\LL(\hat w,v).
	\]
	This gives $\LL(w,v)=\varphi_\ell(v)$.
	By definition of $V(w)$, $\LL(w,v)=p(w)$ follows.
	Hence, we have shown $p(w)=\varphi_\ell(v)$.
	The fact that $v$ solves \eqref{eq:D_NLP_Lagrange} follows
	from assertion~\ref{item:Lagrange_weak}.
	Assertion~\ref{item:Lagrange_saddle_point} is immediate
	from \cite[Theorems~4.7, 4.8, and 4.9]{Ruszczynski2011}.
\end{proof}

\subsubsection{Wolfe duality}\label{sec:Wolfe_duality}

In \cite{Wolfe1961}, Wolfe suggests the dual problem
\begin{equation}\label{eq:D_NLP_Wolfe}\tag{D$_\textup{w}$}
  	\max\limits_{\hat w,\hat v}\{\LL(\hat w,\hat v)\,|\,\nabla_{1}\LL(\hat w,\hat v)=0,\,\hat v\geq 0\}
\end{equation}
associated with \eqref{eq:NLP}.
For feasible points $w\in\R^s$ of \eqref{eq:NLP} and $(\hat w,\hat v)\in\R^s\times\R^t$ of \eqref{eq:D_NLP_Wolfe},
convexity of $p$ and $q_1,\ldots,q_t$ yields
\begin{equation}\label{eq:weak_Wolfe_duality}
\begin{aligned}
	p(w)
	&\geq
	p(\hat w)+p'(\hat w)(w-\hat w)
	=
	p(\hat w)-\hat v^\top q'(\hat w)(w-\hat w)
	\\
	&\geq
	p(\hat w)+\hat v^\top(q(\hat w)-q(w))
	\geq
	p(\hat w)+\hat v^\top q(\hat w)
	=
	\LL(\hat w,\hat v),
\end{aligned}
\end{equation}
i.e., a weak duality relation holds between \eqref{eq:NLP} and \eqref{eq:D_NLP_Wolfe}.
Given validity of a constraint qualification, a strong duality relation can be shown.
The following proposition summarizes the essential duality relations between
\eqref{eq:NLP} and \eqref{eq:D_NLP_Wolfe}.

\begin{proposition}\label{prop:Wolfe_duality}
	\begin{enumerate}
	\item\label{item:Wolfe_weak}
		For feasible points $w\in\R^s$ of \eqref{eq:NLP}
		and $(\hat w,\hat v)\in\R^s\times\R^t$ of \eqref{eq:D_NLP_Wolfe},
		we have $p(w)\geq\LL(\hat w,\hat v)$.
		Particularly, the infimal value of \eqref{eq:NLP} is not smaller
		than the supremal value of \eqref{eq:D_NLP_Wolfe}.
	\item\label{item:Wolfe_strong}
		Let $w\in\R^s$ be a minimizer of \eqref{eq:NLP}
		where GCQ holds.
		Then there exists $v\in\R^t$ such that $(w,v)$
		is a global maximizer of \eqref{eq:D_NLP_Wolfe}
		and $p(w)=\LL(w,v)$.
	\item\label{item:Wolfe_strong_converse}
		Let $(\hat w,\hat v)\in\R^s\times\R^t$ be a global maximizer of \eqref{eq:D_NLP_Wolfe}
		and assume that $\nabla^2_{1,1}\LL(\hat w,\hat v)$ is regular.
		Then $\hat w$ is a minimizer of \eqref{eq:NLP}, and $\hat v\in V(\hat w)$.
	\end{enumerate}
\end{proposition}
\begin{proof}
	Assertion~\ref{item:Wolfe_weak} has been shown above,
	and a proof of assertion~\ref{item:Wolfe_strong_converse}
	can be found in \cite{CravenMond1971}.
	
	In order to show assertion~\ref{item:Wolfe_strong},
	one has to realize that validity of GCQ ensures that $w$ is
	a stationary point of \eqref{eq:NLP}, i.e.,
	there exists $v\in V(w)$.
	Hence, $(w,v)$ is a saddle point of the Lagrangian $\LL$.
	Following the proof of \cite[Theorem~2]{Wolfe1961} now yields the claim.	
\end{proof}

Note that, for a global maximizer $(\hat w,\hat v)\in\R^s\times\R^t$ of \eqref{eq:D_NLP_Wolfe},
$\hat w$ does not need to be a minimizer of \eqref{eq:NLP}, i.e.,
\cref{prop:Wolfe_duality}\,\ref{item:Wolfe_strong_converse} does not generally hold
in the absence of the additional regularity condition.
Clearly, the latter is inherently satisfied whenever the function $p$ is strictly convex
as this yields
positive definiteness of the Hessian $\nabla^2p(w)$ for each $w\in\R^s$,
while the matrices $\nabla^2q_i(w)$, $i=1,\ldots,t$, are positive semidefinite.
A direct proof of \cref{prop:Wolfe_duality}\,\ref{item:Wolfe_strong_converse}
under this stronger assumption is presented in \cite[Theorem~3.3]{Kanniappan1984}.
Observe that one does not need validity of a constraint qualification to distill this result,
as claimed in \cite[Theorem~2.3]{LiLinZhangZhu2022}.

We emphasize that assertions~\ref{item:Wolfe_weak}
and~\ref{item:Wolfe_strong} of \cref{prop:Wolfe_duality}
remain true whenever the mapping $\LL(\cdot,v)$ is so-called pseudoconvex for each $v\in\R^t_+$,
see \cite[Theorems~2.1, 2.2]{LiLinZhangZhu2022}.
This is clearly milder than our inherent convexity assumption
on the functions $p$ and $q_1,\ldots,q_t$.
For simplicity of presentation, however, we will rely on this stronger assumption.

\subsubsection{Mond--Weir duality}\label{sec:Mond_Weir_duality}

In \cite{MondWeir1981}, Mond and Weir suggest the consideration of the dual problem
\begin{equation}\label{eq:D_NLP_MondWeir}\tag{D$_\textup{mw}$}
	\max\limits_{\hat w,\hat v}\{
		p(\hat w)\,|\,\nabla_1\LL(\hat w,\hat v)=0,\,\hat v^\top q(\hat w)\geq 0,\,\hat v\geq 0
		\}
\end{equation}
associated with \eqref{eq:NLP}.
For feasible points $w\in\R^s$ of \eqref{eq:NLP}
and $(\hat w,\hat v)\in\R^s\times\R^t$ of \eqref{eq:D_NLP_MondWeir},
we can show $p(w)\geq p(\hat w)$ as in \eqref{eq:weak_Wolfe_duality}.
One merely has to replace the last equation by the lower estimate ``$\geq p(\hat w)$''
which holds due to $\hat v^\top q(\hat w)\geq 0$.
Hence, we observe a weak duality relation between \eqref{eq:NLP} and \eqref{eq:D_NLP_MondWeir}.
This can also be distilled from \cite[Theorem~3]{EguroMond1986}.
In the subsequent proposition, we summarize essential duality relations between those problems.

\begin{proposition}\label{prop:MondWeir_duality}
	\begin{enumerate}
	\item\label{item:MondWeir_weak}
		For feasible points $w\in\R^s$ of \eqref{eq:NLP}
		and $(\hat w,\hat v)\in\R^s\times\R^t$ of \eqref{eq:D_NLP_MondWeir},
		we have $p(w)\geq p(\hat w)$.
		Particularly, the infimal value of \eqref{eq:NLP} is not smaller
		than the supremal value of \eqref{eq:D_NLP_MondWeir}.
	\item\label{item:MondWeir_strong}
		Let $w\in\R^s$ be a minimizer of \eqref{eq:NLP}
		where GCQ holds.
		Then there exists $v\in\R^t$ such that $(w,v)$
		is a global maximizer of \eqref{eq:D_NLP_MondWeir}.
	\item\label{item:MondWeir_strong_converse}
		Let $(\hat w,\hat v)\in\R^s\times\R^t$
		be a global maximizer of \eqref{eq:D_NLP_MondWeir}
		and assume that
		$\nabla^2_{1,1}\LL(\hat w,\hat v)$ is positive definite.
		Furthermore, assume that one of the following conditions holds.
		\begin{enumerate}
			\item[(i)] The gradient $\nabla p(\hat w)$ does not vanish.
			\item[(ii)] Problem \eqref{eq:NLP} possesses a minimizer
				where GCQ holds.
		\end{enumerate}
		Then $\hat w$ is a minimizer of \eqref{eq:NLP}, and $\hat v\in V(\hat w)$.
	\end{enumerate}
\end{proposition}
\begin{proof}
	Assertion~\ref{item:MondWeir_weak} has been shown above,
	and assertion~\ref{item:MondWeir_strong} follows from \cite[Theorem~4]{EguroMond1986}.
	
	Let us prove assertion~\ref{item:MondWeir_strong_converse} under the additional
	assumption $\nabla p(\hat w)\neq 0$.
	We apply the Fritz--John conditions,
	see \cite[Proposition~3.3.5]{Bertsekas1999}, to \eqref{eq:D_NLP_MondWeir} in order to find
	multipliers $(\alpha,\beta,\gamma,\delta)\in\R\times\R^s\times\R\times\R^t$ such that
	\begin{subequations}\label{eq:FJ_MW}
		\begin{align}
			\label{eq:FJ_MW_w}
				-\alpha\nabla p(\hat w) + \nabla^2_{1,1}\LL(\hat w,\hat v)\beta
				-\gamma q'(\hat w)^\top\hat v
				&=
				0,
				\\
			\label{eq:FJ_MW_v}
				q'(\hat w)\beta-\gamma q(\hat w) - \delta
				&=
				0,
				\\
			\label{eq:FJ_MW_cs_gamma}
				\gamma\geq 0,\quad\gamma\,\hat v^\top q(\hat w)
				&=
				0,
				\\
			\label{eq:FJ_MW_cs_delta}
				\delta\geq 0,\quad\delta^\top\hat v
				&=
				0,
				\\
			\label{eq:FJ_MW_nontrivial}
				|\alpha|+\norm{\beta}+|\gamma|+\norm{\delta}
				&>
				0.
		\end{align}
	\end{subequations}
	Some rearrangements in \eqref{eq:FJ_MW_w} yield
	\begin{align*}
		\nabla^2_{1,1}\LL(\hat w,\hat v)\beta
		&=
		\alpha\nabla p(\hat w)+\gamma q'(\hat w)^\top\hat v
		\\
		&=
		\alpha\nabla_1\LL(\hat w,\hat v) + (\gamma-\alpha) q'(\hat w)^\top\hat v
		=
		(\gamma-\alpha) q'(\hat w)^\top\hat v
	\end{align*}	
	since $(\hat w,\hat v)$ is feasible to \eqref{eq:D_NLP_MondWeir}.
	Multiplying this equation from the left by $\beta^\top$,
	positive definiteness of $\nabla^2_{1,1}\LL(\hat w,\hat v)$
	together with \eqref{eq:FJ_MW_v}, \eqref{eq:FJ_MW_cs_gamma},
	and \eqref{eq:FJ_MW_cs_delta} yields
	\begin{align*}
		0
		&\leq
		\beta^\top\nabla^2_{1,1}\LL(\hat w,\hat v)\beta
		=
		(\gamma-\alpha) (q'(\hat w)\beta)^\top\hat v
		\\
		&=
		(\gamma-\alpha)(\gamma q(\hat w)+\delta)^\top\hat v
		=
		(\gamma-\alpha)\,\gamma\,q(\hat w)^\top \hat v
		+
		(\gamma-\alpha)\,\delta^\top\hat v
		=
		0,
	\end{align*}
	which is why we find $\beta=0$.
	
	Supposing $\gamma=0$, \eqref{eq:FJ_MW_v} yields $\delta=0$,
	and \eqref{eq:FJ_MW_w} then also gives $\alpha=0$ since $\nabla p(\hat w)\neq 0$.
	This, however, is a contradiction to \eqref{eq:FJ_MW_nontrivial}.
	Hence, we have $\gamma>0$, so that \eqref{eq:FJ_MW_v} and \eqref{eq:FJ_MW_cs_delta}
	guarantee $q(\hat w)=-\delta/\gamma\leq 0$.
	Furthermore, \eqref{eq:FJ_MW_cs_gamma} yields $\hat v^\top q(\hat w)=0$.
	As $\nabla_1\LL(\hat w,\hat v)=0$ is also valid,
	$\hat w$ is a feasible and stationary point of \eqref{eq:NLP}.
	Convexity of \eqref{eq:NLP}, thus, implies that $\hat w$ is a minimizer 
	of \eqref{eq:NLP}.
	
	Note that $\nabla_1\LL(\hat w,\hat v)=0$ implies that $\hat w$ is a minimizer
	of the convex function $w\mapsto\LL(w,\hat v)$.
	Hence, positive definiteness of the Hessian $\nabla^2_{1,1}\LL(\hat w,\hat v)$ ensures
	that $\hat w$ is the uniquely determined minimizer of this function.
	Consequently, under the second additional assumption,
	assertion~\ref{item:MondWeir_strong_converse} follows immediately from \cite{GulatiCraven1983}.
\end{proof}

Let us note that \cref{prop:MondWeir_duality}\,\ref{item:MondWeir_strong_converse}
puts the main result in \cite{GulatiCraven1983} as well as \cite[Theorem~7]{EguroMond1986}
and, thus, \cite[Theorem~2.3]{LiLinZhu2024} into some new light as one does not necessarily need
to assume validity of a constraint qualification in the primal problem \eqref{eq:NLP}
in order to infer so-called strict converse duality results in the setting of Mond--Weir duality.
It is enough to assume that either $p$ or one of the active constraint functions $q_i$ associated
with a positive Lagrange multiplier $\hat v_i$ is strictly convex
in order to ensure validity of the definiteness assumption,
and one has to postulate that the gradient of $p$ at the reference point is nonvanishing.
We also want to emphasize that, in the context of our simplified setting,
the proof of \cite[Theorem~7]{EguroMond1986} requires strict convexity of $p$.
The assumption $\nabla p(\hat w)\neq 0$ is not required
in \cite{EguroMond1986,GulatiCraven1983}.
Let us hint the interested reader to the fact that, in \cite{EguroMond1986}, the authors work with so-called
strict pseudo-convexity of the objective function,
but the definition of this property in \cite[Definition~5]{EguroMond1986} is not correct
as one can see from the proof of \cite[Theorem~7]{EguroMond1986}.

It remains open whether one can omit the additional assumptions in
\cref{prop:MondWeir_duality}\,\ref{item:MondWeir_strong_converse}.
Let us also hint the reader to the fact that the assumptions in
\cref{prop:MondWeir_duality}\,\ref{item:MondWeir_strong_converse}
are stronger than those ones used in
\cref{prop:Wolfe_duality}\,\ref{item:Wolfe_strong_converse}.

Let us close the section by remarking that \cref{prop:MondWeir_duality}\,\ref{item:MondWeir_weak}
and~\ref{item:MondWeir_strong} can be shown to hold
in the presence of suitable pseudo- and quasi-convexity properties
of the involved data functions, see \cite[Theorems~3 and~4]{EguroMond1986}.
For simplicity of presentation, however, we restrict ourselves to the fully convex situation
to avoid technicalities.

\section{The optimistic bilevel optimization problem and its standard reformulations}\label{sec:bpp_intro}

For sufficiently smooth functions
$f\colon\R^n\times\R^m\to\R $ and $g\colon\R ^n\times\R ^m\to\R^p$,
we consider the lower-level parametric optimization problem
\begin{equation}\label{eq:llp}\tag{P$(x)$}
  \min\limits_y \{f(x,y)\,|\, g(x,y)\leq 0\},
\end{equation}
its optimal value function $\varphi\colon\R^n\to\overline\R$ given by
\begin{equation*}
  \forall x\in\R^n\colon\quad
  \varphi(x):=\inf\limits_y \{f(x,y)\,|\, g(x,y)\leq 0\},
\end{equation*}
and its solution mapping $\Psi\colon\R^n\tto\R^m$ defined via
\begin{equation*}
  \forall x\in\R^n\colon\quad
  \Psi(x):=\{y\in\R^m\,|\, g(x,y)\leq 0,\,f(x,y)\leq \varphi(x)\}.
\end{equation*}
Then, again for a sufficiently smooth function $F\colon\R^n\times\R^m\to\R$
and some nonempty, closed set $X\subset\R^n$,
the upper-level optimization problem
\begin{equation}\label{eq:ulp}\tag{BOP}
  \textrm{``}\min\limits_x\textrm{''} \{F(x,y)\,|\,x\in X,\ (x,y)\in \gph \Psi\}
\end{equation}
can be defined.
Due to its two-level structure,
we refer to \eqref{eq:ulp} as a bilevel optimization problem.
	Typically, the set $X$ is modeled via inequality constraints.
Let us note that equality constraints can be easily added
to \eqref{eq:llp} and \eqref{eq:ulp} without enhancing the difficulty of the problem.

Problem \eqref{eq:ulp} is not well-defined in the case
where $\Psi$ is not singleton-valued on the intersection of its domain and $X$, i.e.,
when there exists $x\in X$ such that problem \eqref{eq:llp}
has more than one globally optimal solution.
Two ways out are investigated in such a situation,
namely the optimistic and the pessimistic bilevel optimization problem.
In the optimistic bilevel optimization problem
\begin{equation}\label{eq:oblp}
  \min\limits_x\{\inf\limits_y \{F(x,y)\,|\, y\in\Psi(x)\}\,|\,x\in X\},
\end{equation}
the leader assumes cooperation of the follower
in the sense that the follower will select an optimal solution of problem \eqref{eq:llp}
which is a best one w.r.t.\ the leader's objective function.
If this is not the case (or not allowed),
the leader is forced to bound the damage resulting from an unwelcome selection of the follower,
which results in the pessimistic bilevel optimization problem
\begin{equation*}
  \min\limits_x\{\sup\limits_y \{F(x,y)\,|\, y\in\Psi(x)\}\,|\,x\in X\}.
\end{equation*}
Let us note that problem \eqref{eq:oblp} is almost equivalent to
\begin{equation}\label{eq:eoblp}\tag{OBOP}
  \min\limits_{x,y}\{F(x,y)\,|\, x\in X,\, y\in \Psi(x)\},
\end{equation}
see \cite[Proposition~6.9]{DempeMordukhovichZemkoho2012} for a precise analysis.
Throughout the paper, we will focus our attention to the model \eqref{eq:eoblp}
which, in the literature, is often referred to as an optimistic bilevel optimization
problem as well.

Let us mention that, by definition of $\Psi$,
\eqref{eq:eoblp} can be equivalently stated as the single-level optimization problem
\begin{equation}\label{eq:VFR}\tag{VF$_\textup{ref}$}
	\min\limits_{x,y}\{F(x,y)\,|\, x\in X,\,g(x,y)\leq 0,\,f(x,y)\leq\varphi(x)\},
\end{equation}
the so-called value function reformulation of \eqref{eq:eoblp}, see \cite{Outrata1990}.
Due to the appearance of the implicitly known function $\varphi$,
which is likely to be nonsmooth in several practically relevant situations,
\eqref{eq:VFR} is still a rather challenging problem.
Besides, NSMFCQ is likely to be violated at the feasible points of this problem,
see \cite[Proposition~3.2]{YeZhu1995}.

Subsequently, let $L\colon\R^n\times\R^m\times\R^p\to\R$ given by
\[
	\forall x\in\R^n,\,\forall y\in\R^m,\,\forall u\in\R^p\colon\quad
	L(x,y,u):=f(x,y)+u^\top g(x,y)
\]
denote the Lagrangian function associated with the lower-level problem \eqref{eq:llp}.
Furthermore, we make use of the lower-level Lagrange multiplier mapping
$\Lambda\colon\R^n\times\R^m\tto\R^p$ defined via
\[
	\forall x\in\R^n,\,\forall y\in\R^m\colon\quad
	\Lambda(x,y):=\{u\in\R^p\,|\,\nabla_2L(x,y,u)=0,\,0\leq u \perp -g(x,y)\geq 0\}.
\]
By continuous differentiability of the data functions $f$ and $g$,
$\Lambda$ possesses a closed graph.
Whenever the functions $f(x,\cdot),g_1(x,\cdot),\ldots,g_p(x,\cdot)$
are convex for each $x\in X$, which we will assume subsequently,
then the inclusion $\dom\Lambda\subset\gph\Psi$ is generally valid.
The converse is true whenever, for each $(x,y)\in\gph\Psi$,
GCQ holds for \eqref{eq:llp} at $y$.
This observation motivates the investigation of the single-level surrogate problem
\begin{equation}\label{eq:KKTref}\tag{KKT$_\text{ref}$}
	\begin{aligned}
  	&\min\limits_{x,y,u}&	&F(x,y)&		\\
  	&\subjectto&			&x\in X,&		\\
  	&&						&g(x,y)\leq 0,\,u\geq 0,\,u^\top g(x,y)=0,&	\\
  	&&						&\nabla_2L(x,y,u)=0,&	
  \end{aligned}
\end{equation}
which is referred to as the KKT reformulation of \eqref{eq:eoblp}.
We note that the latter problem is a so-called
mathematical problem with complementarity constraints,
and, consequently, MFCQ is violated at each feasible point of \eqref{eq:KKTref},
see e.g.\ \cite[Proposition~1.1]{YeZhuZhu1997}.
Additionally, the lower-level Lagrange multiplier $u$ is 
a variable in \eqref{eq:KKTref} since the implicit problem
\[
	\min\limits_{x,y}\{F(x,y)\,|\,x\in X,\,(x,y)\in\dom\Lambda\}
\]
is seemingly harder to treat.
The presence of the additional variable $u$ in \eqref{eq:KKTref}, however,
causes that \eqref{eq:eoblp} and its transformation \eqref{eq:KKTref} are no
longer equivalent problems, even if a constraint qualification is satisfied
at each feasible point of the lower-level problem \eqref{eq:llp} for each $x\in X$.
In fact, \eqref{eq:KKTref} is likely to possess artificial local minimizers and
stationary points which are not related to local minimizers and stationary points
of \eqref{eq:eoblp}. This phenomenon has been reported in the classical paper
\cite{DempeDutta2012}. Furthermore, in general, qualification conditions which have to
be employed to tackle \eqref{eq:KKTref} in order to derive necessary optimality conditions
are, in general, stronger than those ones needed
to access \eqref{eq:eoblp} via its fully equivalent transformation
\begin{equation}\label{eq:GEref}\tag{GE$_\textup{ref}$}
	\min\limits_{x,y}\{F(x,y)\,|\,x\in X,\,-\nabla_2f(x,y)\in N_{\Gamma(x)}(y)\},
\end{equation}
see \cite{AdamHenrionOutrata2018},
where we used the lower-level feasibility mapping $\Gamma\colon\R^n\tto\R^m$ given by
\begin{equation}\label{eq:feasibility_mapping}
	\forall x\in\R^n\colon\quad
	\Gamma(x):=\{y\in\R^m\,|\,g(x,y)\leq 0\}.
\end{equation}
In \eqref{eq:GEref}, which is called the generalized equation reformulation of \eqref{eq:eoblp},
the variable $u$ is implicitly present as we have
\[
	N_{\Gamma(x)}(y)
	=
	\{ g'_2(x,y)^\top u\in\R^m\,|\,u\geq 0,\,u^\top g(x,y) = 0\}
\]
for each $(x,y)\in(X\times\R^m)\cap\gph\Gamma$ such that GCQ holds for \eqref{eq:llp} at $y$.
That is why $u$ is referred to as an implicit variable in \eqref{eq:GEref}.
Making it explicit, one ends up with problem \eqref{eq:KKTref} which is seemingly easier to
treat as one can omit dealing with generalized derivatives of the normal cone mapping
$(x,y)\mapsto N_{\Gamma(x)}(y)$. The latter, however, is nowadays comparatively well understood,
see e.g.\ \cite{GfrererOutrata2016a,GfrererOutrata2016b}, and these works also influenced the
aforementioned paper \cite{AdamHenrionOutrata2018}.
Recently, the phenomenon of implicit variables has been embedded into a much more general
framework in \cite{BenkoMehlitz2021} where it is rigorously analyzed from the perspective
of limiting variational analysis.

\section{Duality-based transformations of the optimistic bilevel optimization problem}\label{sec:duality_based_trafos}

In this section, we introduce and study three different single-level reformulations
of the bilevel optimization problem \eqref{eq:eoblp} which are constructed by replacing
the implicit constraint $(x,y)\in\gph\Psi$ by enforcing strong duality for the
lower-level problem in the sense of Lagrange, Wolfe, and Mond--Weir
in \cref{sec:trafo_Lagrange_duality,sec:trafo_Wolfe_duality,sec:trafo_Mond_Weir_duality},
respectively.
More precisely, in the reformulated problems, we state constraints which ensure
feasibility for the lower-level problem and the associated dual as well as
validity of strong duality. These transformations come for the price of additional
variables, namely, the variables of the lower-level dual problem.
\cref{sec:comparison} briefly summarizes the obtained results 
by presenting a comparison.

Throughout the section, we assume that the data functions
$f(x,\cdot),g_1(x,\cdot),\ldots,g_p(x,\cdot)$ are convex for each $x\in X$.
Hence, \eqref{eq:llp} is a convex optimization problem for each $x\in X$,
and we can access it via the duality theory from \cref{sec:duality}. 
As already mentioned in \cref{sec:duality}, the necessary duality relations may also
hold in more general settings where the lower-level data functions are, in a certain sense,
pseudo- or quasi-convex in the follower's variables, but for simplicity of presentation,
we focus on the standard convex situation. The latter is sufficient to illustrate all
features of the duality-based single-level reformulations of \eqref{eq:eoblp}
under consideration.
We emphasize that none of the three discussed reformulations is reasonable as soon as
the lower-level problem and its dual do not enjoy at least a weak duality relation.
For simplicity of presentation,
we also assume that $f$ and $g$ are twice continuously differentiable
although for the majority of our results, it is enough that $f$ and $g$ are just
continuously differentiable.

\subsection{Transformation based on Lagrange duality}\label{sec:trafo_Lagrange_duality}

We start by considering a single-level reformulation which is based on the lower-level
Lagrange dual problem, see \cref{sec:Lagrange_duality}.
Therefore, let us introduce the Lagrange value function
$\psi_\ell\colon\R^n\times\R^p\to\overline{\R}$ of \eqref{eq:llp} by means of
\[
	\forall x\in\R^n,\,\forall u\in\R^p\colon\quad
	\psi_\ell(x,u)
	:=
	\inf\limits_y\{L(x,y,u)-\delta_{\R^p_+}(u)\,|\,y\in\R^m\}.
\]
Similarly as in \cref{sec:Lagrange_duality}, one can show that this function is upper semicontinuous.
Now, a potential single-level reformulation of \eqref{eq:eoblp} is given by
\begin{equation}\label{eq:LDref}\tag{LD$_\textup{ref}$}
  \begin{aligned}
  	&\min\limits_{x,y,u}&	&F(x,y)&		\\
  	&\subjectto&			&x\in X,&		\\
  	&&						&g(x,y)\leq 0,\,u\geq 0,&	\\
  	&&						&f(x,y)\leq \psi_\ell(x,u),&	
  \end{aligned}
\end{equation}
and we will refer to the latter as the Lagrange dual reformulation of \eqref{eq:eoblp}.
By upper semicontinuity of $\psi_\ell$, the feasible set of \eqref{eq:LDref} is closed.
This reformulation approach has been suggested in \cite{OuattaraAswani2018}
where the authors, additionally, assume that the lower-level feasible set $\Gamma(x)$
is uniformly bounded for each $x\in X$. Here, we do not postulate such a restrictive assumption.
Note that, on the one hand, for each feasible point $(\bar x,\bar y,\bar u)\in\R^n\times\R^m\times\R^p$
of \eqref{eq:LDref}, \cref{prop:Lagrange_duality}\,\ref{item:Lagrange_weak} ensures
that $\bar y$ solves \hyperref[eq:llp]{\textup{(P$(\bar x)$)}}
to optimality while $\bar u$ is a maximizer of
\[
	\max\limits_{u}\{\psi_\ell(\bar x,u)\,|\,u\geq 0\}.
\]
Particularly, $(\bar x,\bar y)$ is feasible to \eqref{eq:eoblp}.
On the other hand, since \cref{prop:Lagrange_duality}\,\ref{item:Lagrange_weak}
guarantees $f(x,y)\geq\psi_\ell(x,u)$ for each feasible triplet
$(x,y,u)\in\R^n\times\R^m\times\R^p$ of \eqref{eq:LDref}, taking the infimum over
all $y\in\Gamma(x)$ yields $\varphi(x)\geq\psi_\ell(x,u)$.
Furthermore, $\varphi(x)\leq f(x,y)\leq\psi_\ell(x,u)$ is immediately clear
by definition of the lower-level optimal value function $\varphi$.
Hence, $\varphi(x)=\psi_\ell(x,u)$ follows,
and this shows that \eqref{eq:LDref} is closely related to \eqref{eq:VFR}.
More specifically, $(x,y)$ is feasible to \eqref{eq:VFR}.
Conversely, there may exist feasible points of \eqref{eq:VFR} which do not correspond
to feasible points of \eqref{eq:LDref} as soon as strong Lagrange duality fails to hold
for the lower-level problem.
We also want to note that, by definition, the Lagrange value function $\psi_\ell$
is inherently discontinuous while there is a fairly good chance
for the lower-level optimal value function $\varphi$ to be locally
Lipschitz continuous under comparatively mild assumptions,
see e.g.\ \cite[Corollary~4.8]{GuoLinYeZhang2014} or \cite[Section~4.1.1]{MehlitzMinchenko2022}
where local Lipschitzness of $\varphi$ is proven in the presence of the
quasi-normality and relaxed constant positive linear dependence constraint qualification,
respectively.

We aim to study the relationship between \eqref{eq:eoblp} and its reformulation
w.r.t.\ global and local minimizers. Noting that the additional variable $u$ in \eqref{eq:LDref}
plays the role of an implicit variable which has been made explicit,
the situation at hand is a particular
setting of the far more general framework studied in \cite{BenkoMehlitz2021}.
In order to apply it here, we introduce a so-called intermediate mapping
$K_\ell\colon\R^n\times\R^m\tto\R^p$ by means of
\[
	\forall x\in\R^n,\,\forall y\in\R^m\colon\quad
	K_\ell(x,y)
	:=
	\{u\in\R^p\,|\,g(x,y)\leq 0,\,u\geq 0,\,f(x,y)\leq\psi_\ell(x,u)\}.
\]
In the lemma below, we study the properties of this mapping.
\begin{lemma}\label{lem:properties_intermediate_map_Lagrange}
	\begin{enumerate}
		\item\label{item:K_l_closed} The set $\gph K_\ell$ is closed.
		\item\label{item:K_l_Lagrange_multipliers}
			For each pair $(x,y)\in\R^n\times\R^m$,
			we have $K_\ell(x,y)=\Lambda(x,y)$.
		\item\label{item:K_l_dom} Fix $(\bar x,\bar y)\in\gph\Psi$ and assume that
			GCQ
			holds for \hyperref[eq:llp]{\textup{(P$(\bar x)$)}} at $\bar y$.
			Then $(\bar x,\bar y)\in\dom K_\ell$.
	\end{enumerate}
\end{lemma}
\begin{proof}
	Assertion~\ref{item:K_l_closed} is an immediate consequence of the upper semicontinuity of $\psi_\ell$,
	and assertion~\ref{item:K_l_Lagrange_multipliers} follows from
	\cref{prop:Lagrange_duality}\,\ref{item:Lagrange_weak} and~\ref{item:Lagrange_saddle_point}.
	Finally, assertion~\ref{item:K_l_dom} is a consequence of assertion~\ref{item:K_l_Lagrange_multipliers}.
\end{proof}

Let us hint the reader to the fact that,
due to \cref{lem:properties_intermediate_map_Lagrange}\,\ref{item:K_l_Lagrange_multipliers},
\eqref{eq:LDref} and \eqref{eq:KKTref} possess the same feasible sets.
On the one hand, this means that, comparing \eqref{eq:eoblp} and \eqref{eq:LDref} w.r.t.\ their global
and local minimizers, we will end up with more or less the same results which have been stated in
\cite[Sections~2 and 3]{DempeDutta2012}.
On the other hand, the different formulation of the feasible set in \eqref{eq:LDref} compared to the one
of \eqref{eq:KKTref} may yield better chances for the validity of a constraint qualification
or may open the door to novel algorithmic approaches for the numerical solution of \eqref{eq:eoblp}.

To start, let us point out the relationship between the global minimizers of \eqref{eq:eoblp} and \eqref{eq:LDref}.
Similar findings can be found in \cite[Theorems~2.1 and 2.3]{DempeDutta2012} and
\cite[Proposition~4]{OuattaraAswani2018}.
In a more general setting, such a result has been presented in \cite[Theorem~4.3]{BenkoMehlitz2021}.

\begin{theorem}\label{thm:LDref_global}
	\begin{enumerate}
		\item\label{item:LDref_global}
			Let $(\bar x,\bar y)\in\R^n\times\R^m$ be a global minimizer of \eqref{eq:eoblp}.
			Then, for each $\bar u\in K_\ell(\bar x,\bar y)$,
			$(\bar x,\bar y,\bar u)$ is a global minimizer of \eqref{eq:LDref}.
		\item Let $(\bar x,\bar y,\bar u)\in\R^n\times\R^m\times\R^p$
			be a global minimizer of \eqref{eq:LDref}.
			Assume that, for each $x\in X\cap\dom\Psi$,
			GCQ holds for \eqref{eq:llp}
			at all points in $\Psi(x)$.
			Then $(\bar x,\bar y)$ is a global minimizer of \eqref{eq:eoblp}.
	\end{enumerate}
\end{theorem}
\begin{proof}
	We start with the proof of the first assertion.
	Let us assume that there is some $\bar u\in K_\ell(\bar x,\bar y)$
	such that $(\bar x,\bar y,\bar u)$ is not a global minimizer of
	\eqref{eq:LDref}. Then there exists a feasible point
	$(\tilde x,\tilde y,\tilde u)\in\R^n\times\R^m\times\R^p$
	of \eqref{eq:LDref} such that $F(\bar x,\bar y)>F(\tilde x,\tilde y)$.
	As \cref{prop:Lagrange_duality}\,\ref{item:Lagrange_weak} yields $\tilde y\in\Psi(\tilde x)$,
	$(\tilde x,\tilde y)$ is a feasible point of \eqref{eq:eoblp},
	which means that $(\bar x,\bar y)$ cannot be a global minimizer of \eqref{eq:eoblp}.
	
	We now proceed with the proof of the second assertion.
	Due to \cref{prop:Lagrange_duality}\,\ref{item:Lagrange_weak},
	$(\bar x,\bar y)$ is feasible to \eqref{eq:eoblp}.
	If $(\bar x,\bar y)$ is not a global minimizer of \eqref{eq:eoblp},
	there exists some feasible point $(\tilde x,\tilde y)\in\R^n\times\R^m$ of \eqref{eq:eoblp}
	such that $F(\bar x,\bar y)>F(\tilde x,\tilde y)$.
	Due to $\tilde x\in X\cap\dom\Psi$ and the postulated assumption,
	\cref{lem:properties_intermediate_map_Lagrange}\,\ref{item:K_l_dom} guarantees the existence of
	$\tilde u\in K_\ell(\tilde x,\tilde y)$.
	Hence, $(\tilde x,\tilde y,\tilde u)$ is feasible to \eqref{eq:LDref}
	which shows that $(\bar x,\bar y,\bar u)$ cannot be a global minimizer
	of that problem.
\end{proof}

We would like to hint the reader to the fact that we do not necessarily have
$\gph\Psi=\dom K_\ell$ in the absence of a constraint qualification
which is why we cannot rely on the proof of \cite[Theorem~4.3]{BenkoMehlitz2021}
as a justification of \cref{thm:LDref_global}.
In fact, if, for some pair $(\bar x,\bar y)\in\gph\Psi$,
\hyperref[eq:llp]{\textup{(P$(\bar x)$)}} does not satisfy a constraint qualification
at $\bar y$, then $K_\ell(\bar x,\bar y)$ can be empty.
In this case, \cref{thm:LDref_global}\,\ref{item:LDref_global} is not very helpful.
It is, thus, desirable that a constraint qualification is satisfied at all feasible points
of the lower-level problem \eqref{eq:llp} for each $x\in X$ to make the transformation
meaningful.

Let us now turn our attention to the relationship between the local minimizers
of \eqref{eq:eoblp} and \eqref{eq:LDref}.
The following result generalizes \cite[Theorem~3.2]{DempeDutta2012}
as well as \cite[Proposition~5]{OuattaraAswani2018},
and has been motivated by \cite[Theorem~4.5]{BenkoMehlitz2021}.
\begin{theorem}\label{thm:LDref_local}
	\begin{enumerate}
		\item\label{item:LDref_local_simple} Let $(\bar x,\bar y)\in\R^n\times\R^m$ be a local minimizer of \eqref{eq:eoblp}.
			Then, for each $\bar u\in K_\ell(\bar x,\bar y)$,
			$(\bar x,\bar y,\bar u)$ is a local minimizer of \eqref{eq:LDref}.
		\item\label{item:LDref_local}
			Let $(\bar x,\bar y, u)\in\R^n\times\R^m\times\R^p$
			be a local minimizer of \eqref{eq:LDref} for each $u\in K_\ell(\bar x,\bar y)$.
			Assume that, for each $x\in X\cap \dom\Psi$ in a neighborhood of $\bar x$,
			GCQ holds for \eqref{eq:llp}
			at all points in $\Psi(x)$.
			Finally, let $K_\ell$ be inner semicompact at $(\bar x,\bar y)$ w.r.t.\ its domain.
			Then $(\bar x,\bar y)$ is a local minimizer of \eqref{eq:eoblp}.
	\end{enumerate}
\end{theorem}
\begin{proof}
	The proof of the first assertion follows in similar way as the 
	validation of the first statement of \cref{thm:LDref_global}.
	We, thus, only prove the second statement.
	Due to \cref{prop:Lagrange_duality}\,\ref{item:Lagrange_weak},
	$(\bar x,\bar y)$ is feasible to \eqref{eq:eoblp}.
	Suppose that $(\bar x,\bar y)$ is not a local minimizer of \eqref{eq:eoblp}.
	Then we find a sequence $\{(x^k,y^k)\}_{k\in\N}\subset (X\times\R^m)\cap\gph\Psi$
	such that $x^k\to\bar x$, $y^k\to\bar y$, and $F(x^k,y^k)<F(\bar x,\bar y)$ for all $k\in\N$.
	For sufficiently large $k\in\N$, \cref{lem:properties_intermediate_map_Lagrange}\,\ref{item:K_l_dom}
	ensures $(x^k,y^k)\in\dom K_\ell$.
	As $K_\ell$ is assumed to be inner semicompact w.r.t.\ its domain at $(\bar x,\bar y)$,
	we can find a sequence $\{u^k\}_{k\in\N}\subset\R^p$ such that
	$u^k\in K_\ell(x^k,y^k)$ holds for all large enough $k\in\N$
	while $\{u^k\}_{k\in\N}$ is bounded along a subsequence.
	Hence, along a subsequence (without relabeling), we may assume $u^k\to\bar u$
	for some $\bar u\in\R^p$.
	Applying \cref{lem:properties_intermediate_map_Lagrange}\,\ref{item:K_l_closed},
	$\bar u\in K_\ell(\bar x,\bar y)$ follows, i.e.,
	$(\bar x,\bar y,\bar u)$ is feasible to \eqref{eq:LDref}.
	Since all points along the tail of the sequence $\{(x^k,y^k,u^k)\}_{k\in\N}$ possess a smaller
	objective value than $(\bar x,\bar y,\bar u)$,
	the latter cannot be a local minimizer of \eqref{eq:LDref}.	
	This, however, contradicts the theorem's assumptions.
\end{proof}

In order to ensure validity of the inner semicompactness requirement
in \cref{thm:LDref_local}\,\ref{item:LDref_local}, different assumptions can be employed.
In \cite[Theorem~3.2]{DempeDutta2012}, the authors make use of Slater's constraint qualification
which already implies that the images of $\Lambda$ and, thus, of $K_\ell$
are uniformly bounded in a neighborhood of $(\bar x,\bar y)$,
see e.g.\ \cite[Proposition~4.43]{BonnansShapiro2000} or \cite[proof of Theorem~4.2]{Dempe2002}.
In \cite[Proposition~5]{OuattaraAswani2018}, validity of the linear independence constraint
qualification is assumed for that purpose, yielding uniqueness of the multiplier and
continuity of the locally singleton-valued multiplier mapping $\Lambda$, see e.g.\
\cite[Lemma~4.44, Proposition~4.47]{BonnansShapiro2000}.
However, inner semicompactness on its own is much weaker than such local boundedness or even
continuity of $K_\ell$.
In \cite[Proposition~3.2, Theorem~3.3]{GfrererMordukhovich2017},
the authors show that a parametric
version of the constant rank constraint qualification
(which is independent of Slater's constraint qualification)
or so-called Robinson stability
(which is, in general, weaker than Slater's constraint qualification)
are also sufficient for the required inner semicompactness.

It has been visualized in \cite[Example~3.4]{DempeDutta2012}
that local optimality of $(\bar x,\bar y,u)$
for \eqref{eq:KKTref} indeed has to established for all $u\in \Lambda(\bar x,\bar y)$
in order to infer local optimality of $(\bar x,\bar y)$ for \eqref{eq:eoblp}.
Hence, due to \cref{lem:properties_intermediate_map_Lagrange}\,\ref{item:K_l_Lagrange_multipliers},
it is evident that one cannot abstain from postulating local optimality of
$(\bar x,\bar y,u)$ for \eqref{eq:LDref} for all $u\in K_\ell(\bar x,\bar y)$
in \cref{thm:LDref_local}\,\ref{item:LDref_local} in order to guarantee
local optimality of $(\bar x,\bar y)$ for \eqref{eq:eoblp}.
From an algorithmic point of view
this is a huge drawback whenever $K_\ell(\bar x,\bar y)$ is not just a singleton,
as one has to verify local optimality of, potentially,
infinitely many points for the surrogate problem.
For the purpose of completeness, and as a preparation for the upcoming subsections,
let us revisit \cite[Example~3.4]{DempeDutta2012}
to illustrate this behavior.

\begin{example}\label{ex:artificial_local_solution_l}
	Consider the linear lower-level problem
	\begin{equation}\label{eq:llp_ex_lagrange}
    	\min_{y}\{-y\,|\, x+y\leq 1,\, -x+y\leq 1\}
	\end{equation}
	and the upper-level problem
	\begin{equation}\label{eq:ulp_ex_lagrange}
	\min\limits_{x,y}\{(x-1)^2+(y-1)^2\,|\, (x,y)\in \gph\Psi\}.
	\end{equation}
	This problem has the uniquely determined global minimizer
	$(\hat x,\hat y):=(1/2,1/2)$,
	and there do not exist local minimizers different from $(\hat x,\hat y)$.
	
	Let us derive an explicit representation of the associated Lagrange dual reformulation
	\eqref{eq:LDref}. Therefore, we observe that, for each $x\in\R$ and $u\in\R^2_+$, one has
	\begin{align*}
		\psi_\ell(x,u)
		&=
		\inf\limits_y\{-y+u_1(x+y-1)+u_2(-x+y-1)\,|\,y\in\R\}
		\\
		&=
		\begin{cases}
			x(u_1-u_2)-1	&	u_1+u_2=1,
			\\
			-\infty			&	\text{otherwise.}
		\end{cases}
	\end{align*}
	This allows to state \eqref{eq:LDref} in the form
	\begin{equation}\label{eq:LDref_example}
		\begin{aligned}
  		&\min\limits_{x,y,u}&	&(x-1)^2+(y-1)^2&		\\
  		&\subjectto&			&x+y\leq 1,\,-x+y\leq 1,\,u\geq 0,&	\\
  		&&						&-y\leq x(u_1-u_2)-1,\,u_1+u_2=1.&	
  		\end{aligned}
	\end{equation}
	Let us consider the point $(\bar x,\bar y,\bar u):=(0,1,(0,1))$ which is feasible to \eqref{eq:LDref_example}.
	Furthermore, we fix a feasible point $(x,y,u)\in\R\times\R\times\R^2$ of \eqref{eq:LDref_example}
	which is close to $(\bar x,\bar y,\bar u)$.
	The first two constraints in \eqref{eq:LDref_example} give $y\leq 1$, i.e.,
	\[
		1\leq x(u_1-u_2)+y\leq x(u_1-u_2)+1
	\]
	which yields $x(u_1-u_2)\geq 0$.
	As we have $u_1-u_2<0$ for $u$ close enough to $\bar u$, $x\leq 0$ follows.
	Then, however, the objective value of $(x,y,u)$ cannot be smaller than the one of $(\bar x,\bar y,\bar u)$,
	i.e., the latter is a local minimizer of \eqref{eq:LDref_example}.
	However, let us emphasize that $(\bar x,\bar y)$ is not a local minimizer of \eqref{eq:ulp_ex_lagrange}.
	
	One can easily check that $K_\ell(\bar x,\bar y)=\{u\in\R^2_+\,|\,u_1+u_2=1\}$,
	so let us fix $\tilde u:=(1,0)\in K_\ell(\bar x,\bar y)$.
	For each $\varepsilon>0$, $(\bar x+\varepsilon,\bar y-\varepsilon,\tilde u)$
	is feasible to \eqref{eq:LDref_example} and, for each $\varepsilon\in(0,1)$,
	possesses a smaller objective value than $(\bar x,\bar y,\tilde u)$.
	Hence, $(\bar x,\bar y,\tilde u)$ is not a local minimizer of \eqref{eq:LDref_example}.
	This observation underlines that one cannot abstain from postulating local optimality
	for all $u\in K_\ell(\bar x,\bar y)$ in \cref{thm:LDref_local}\,\ref{item:LDref_local}.
	It should be noted that all the other assumptions stated in \cref{thm:LDref_local}\,\ref{item:LDref_local}
	are valid. Indeed, the lower-level problem is linear which guarantees validity of GCQ at each
	feasible point of \eqref{eq:llp_ex_lagrange}.
	Furthermore, we have $K_\ell(x,y)\subset K_\ell(\bar x,\bar u)$ for each pair $(x,y)\in\R\times\R$,
	and this uniform boundedness of $K_\ell$ guarantees inner semicompactness 
	w.r.t.\ its domain at each point of its domain.
\end{example}

In the subsequent theorem,
we point out that NSMFCQ
is likely to fail at the feasible points of \eqref{eq:LDref}.
The proof of this result has been inspired by the one of \cite[Proposition~3.2]{YeZhu1995},
where the authors show that NSMFCQ is violated at the feasible points of the
value function reformulation \eqref{eq:VFR} of \eqref{eq:eoblp} under mild assumptions.

\begin{theorem}\label{thm:violation_MFCQ_Lagrange}
 	Let $(\bar x,\bar y,\bar u)\in\R^n\times\R^m\times\R^p$
 	be a feasible point of \eqref{eq:LDref} such that
 	the qualification condition
 	\begin{equation}\label{eq:BCQ}
  		\partial^\infty \chi(\bar x,\bar y,\bar u)
 		\cap
 		\bigl(-N_{\mathcal F}(\bar x,\bar y,\bar u)\bigr)
 		=
 		\{(0,0,0)\}
 	\end{equation}
 	is valid, where $\chi\colon\R^n\times\R^m\times\R^p\to\overline{\R}$ is the function given by
 	\[	
 		\forall x\in\R^n,\,\forall y\in\R^m,\,\forall u\in\R^p\colon\quad
 		\chi(x,y,u):=f(x,y)-\psi_\ell(x,u),
 	\]
 	and $\mathcal F\subset\R^n\times\R^m\times\R^p$ is the set defined via
 	\[
 		\mathcal F:=\{(x,y,u)\in\R^n\times\R^m\times\R^p\,|\,
 						g(x,y)\leq 0,\,u\geq 0
 					\}.
 	\]
 	Furthermore, assume that $X:=\{x\in\R^n\,|\,G(x)\leq 0\}$ holds
 	for a continuously differentiable function $G\colon\R^n\to\R^q$.
 	Then NSMFCQ is violated for \eqref{eq:LDref} at $(\bar x,\bar y,\bar u)$.
\end{theorem}
\begin{proof}
	Whenever MFCQ fails for the constraint system
	\begin{equation}\label{eq:violation_MFCQ_Lagrange_constraint_system}
		g(x,y)\leq 0,\,u\geq 0
	\end{equation}
	at $(\bar x,\bar y,\bar u)$, then NSMFCQ fails for \eqref{eq:LDref} at $(\bar x,\bar y,\bar u)$ as well.
	Hence, let us assume that MFCQ holds for the above constraint system
	\eqref{eq:violation_MFCQ_Lagrange_constraint_system} at $(\bar x,\bar y,\bar u)$.
	Due to \cref{prop:Lagrange_duality}\,\ref{item:Lagrange_weak}, for each triplet
	$(x,y,u)\in\mathcal F$, we have $\chi(x,y,u)\geq 0$.
	Hence, since $(\bar x,\bar y,\bar u)$ is feasible to \eqref{eq:LDref},
	$\chi(\bar x,\bar y,\bar u)=0$ follows.
	That is why $(\bar x,\bar y,\bar u)$ is a global minimizer of the nonsmooth optimization
	problem
	\[
		\min\limits_{x,y,u}\{\chi(x,y,u)\,|\,(x,y,u)\in\mathcal F\}.
	\]
	We apply \cite[Theorem~6.1(ii)]{Mordukhovich2018}, which is possible due to the presumed
	qualification condition \eqref{eq:BCQ}, and obtain
	\[
		(0,0,0)\in\partial\chi(\bar x,\bar y,\bar u)+N_{\mathcal F}(\bar x,\bar y,\bar u).
	\]
	Since MFCQ holds for the constraint system \eqref{eq:violation_MFCQ_Lagrange_constraint_system}
	at $(\bar x,\bar y,\bar u)$, \cite[Theorem~6.14]{RockafellarWets1998} yields the existence
	of multipliers $\alpha,\beta\in\R^p$ such that
	\begin{align*}
		-(g'_1(\bar x,\bar y)^\top\alpha,g'_2(\bar x,\bar y)^\top\alpha,0)
		-(0,0,-\beta)
		\in
		\partial\chi(\bar x,\bar y,\bar u),
		\\
		\alpha\geq 0,\,\alpha^\top g(\bar x,\bar y)=0,
		\\
		\beta\geq 0,\,\beta^\top\bar u=0.
	\end{align*}
	Equivalently, NSMFCQ fails to hold for \eqref{eq:LDref}
	at $(\bar x,\bar y,\bar u)$.
\end{proof}

Let us note that a relaxation approach can be applied to \eqref{eq:LDref} in order to
enforce validity of NSMFCQ at the feasible points of the associated
single-level transformation under additional assumptions, see \cite[Theorem~3]{OuattaraAswani2018}.

Subsequently, we aim to briefly discuss the qualification condition \eqref{eq:BCQ} in more detail.
To start, let us note that the normal cone to $\mathcal F$ at $(\bar x,\bar y,\bar u)$
can be estimated from above in terms of initial problem data in the presence of a mild constraint qualification
like MFCQ for the constraint system in \eqref{eq:violation_MFCQ_Lagrange_constraint_system},
see \cite[Theorem~6.14]{RockafellarWets1998} again.
Hence, it remains to handle the computation or, at least,
estimation of $\partial^\infty\chi(\bar x,\bar y,\bar u)$.
Of course, the latter reduces to $\{(0,0,0)\}$ whenever $\chi$ is locally Lipschitz continuous at
$(\bar x,\bar y,\bar u)$, and this property holds if $\bar u$ is componentwise positive
and at least one of the functions
$f(\bar x,\cdot),g_1(\bar x,\cdot),\ldots,g_p(\bar x,\cdot)$ is strictly convex
as this yields positive definiteness of the Hessian $\nabla^2_{2,2}L(\bar x,\bar y,\bar u)$
and, hence, the parametric optimization problem
\[
	\min\limits_{y}\{L(x,y,u)-\delta_{\R^p_+}(u)\,|\,y\in\R^m\}
\]
possesses, in a neighborhood of $(\bar x,\bar u)$, a single-valued smooth solution mapping
by the implicit function theorem.
In the case where $\bar u$ possesses a vanishing component, these arguments are no longer suitable.
	However, in the presence of a comparatively mild level-boundedness assumption,
	it is possible to show that the mapping $(x,u)\mapsto\inf_y\{L(x,y,u)\,|\,y\in\R^m\}$
	is locally Lipschitz continuous at $(\bar x,\bar u)$,
	see e.g.\ \cite[Corollary~10.14]{RockafellarWets1998}.
	Then, on the one hand, \cite[Exercise~10.10]{RockafellarWets1998} yields
	\[
		\partial^\infty\chi(\bar x,\bar y,\bar u)
		\subset
		\{(0,0)\}\times\partial^\infty \delta_{\R^p_+}(\bar u)
		=
		\{(0,0)\}\times N_{\R^p_+}(\bar u).
	\]
	On the other hand, 
	\cite[Proposition~6.41]{RockafellarWets1998} can be applied to find
	\begin{equation}\label{eq:normal_cone_to_F}
		N_{\mathcal F}(\bar x,\bar y,\bar u)
		=
		N_{\gph\Gamma}(\bar x,\bar y)\times N_{\R^p_+}(\bar u),
	\end{equation}
	and due to $N_{\R^p_+}(\bar u)\subset\R^p_-$,
	validity of \eqref{eq:BCQ} follows.
	The latter, thus, is a comparatively weak condition.
	Let us sketch yet another way to verify validity of \eqref{eq:BCQ}.
We observe that
\[
	\forall x\in\R^n,\,\forall y\in\R^m,\,\forall u\in\R^p\colon\quad
	\chi(x,y,u)
	=
	f(x,y)+\upsilon(x,u)+\delta_{\R^p_+}(u)
\]
holds for the function $\upsilon\colon\R^n\times\R^p\to\overline{\R}$ which is given by
\[
	\forall x\in\R^n,\,\forall u\in\R^p\colon\quad
	\upsilon(x,u):=\sup\limits_y\{-L(x,y,u)\,|\,y\in\R^m\},
\]
so one could rely on suitable marginal function rules for supremum functions,
see e.g.\ \cite{CorreaHantouteLopez2021,GuoYeZhang2024,HantouteLopezCerda2023,PerezAros2019},
and the sum rule for the
singular subdifferential, see \cite[Exercise~10.10]{RockafellarWets1998} again,
in order to estimate $\partial^\infty\chi(\bar x,\bar y,\bar u)$ from above.
A direct application of subdifferential calculus rules to the original definition
of $\chi$ leads to the issue of estimating $\partial^\infty(-\psi_\ell)(\bar x,\bar u)$
from above. Due to the non-Lipschitzness of $\psi_\ell$ and the fact that the
singular subdifferential does not enjoy a plus-minus symmetry property,
this is a rather involved problem in variational analysis, and we abstain from any further
discussion here in order to avoid distorting the focus of the paper.

 	In the upcoming example, we verify that \eqref{eq:BCQ} even holds in 
	\cref{ex:artificial_local_solution_l} where the aforementioned explicit
	criteria for its validity cannot be applied as the lower-level problem is linear.
	\begin{example}\label{ex:BCQ}
		Let us consider the bilevel optimization problem
		\eqref{eq:ulp_ex_lagrange} from \cref{ex:artificial_local_solution_l}.		
		We already verified that
		\[
			\forall x\in\R,\,\forall u\in\R^2\colon\quad
			\psi_\ell(x,u)
			=
			\begin{cases}
				x(u_1-u_2)-1	&	u\in\Delta^2,\\
				-\infty			&	\text{otherwise}
			\end{cases}
		\]
		holds for the associated lower-level-problem \eqref{eq:llp_ex_lagrange}, 
		where $\Delta^2:=\{u\in\R^2_+\,|\,u_1+u_2=1\}$ is the standard simplex in $\R^2$.
		Hence, we find
		\[
			\forall x\in\R,\,\forall y\in\R,\,\forall u\in\R^2\colon\quad
			\chi(x,y,u)
			=
			\begin{cases}
				-y-x(u_1-u_2)+1	&	u\in\Delta^2,\\
				\infty			&	\text{otherwise.}
			\end{cases}
		\]
		Next, we fix an arbitrary feasible point $(\bar x,\bar y,\bar u)\in\R\times\R\times\R^2$
		of the associated problem \eqref{eq:LDref} which is stated in \eqref{eq:LDref_example}.
		Then we have $\bar u\in\Delta^2$,
		and \cite[Exercise~10.10]{RockafellarWets1998} yields
		\[
			\partial^\infty\chi(\bar x,\bar y,\bar u)
			=
			\{(0,0)\}\times N_{\Delta^2}(\bar u).
		\]
		Taking \eqref{eq:normal_cone_to_F} into account, \eqref{eq:BCQ} reduces to
		\[
			N_{\Delta^2}(\bar u)\cap\bigl(-N_{\R^2_+}(\bar u)\bigr)=\{0\}.
		\]
		As we have
		\begin{align*}
			N_{\Delta^2}(\bar u)
			&=
			\begin{cases}
				\{\eta\in\R^2\,|\,-\eta_1+\eta_2\leq 0\}	&	\bar u=(1,0),
				\\
				\{\eta\in\R^2\,|\,\eta_1=\eta_2\}			&	0<\bar u_1,\bar u_2<1,
				\\
				\{\eta\in\R^2\,|\,\eta_1-\eta_2\leq 0\}		&	\bar u=(0,1),
			\end{cases}
			\\
			N_{\R^2_+}(\bar u)
			&=
			\begin{cases}
				\{0\}\times\R_-				&	\bar u=(1,0),
				\\
				\{(0,0)\}					&	0<\bar u_1,\bar u_2<1,
				\\
				\R_-\times\{0\}				&	\bar u=(0,1),
			\end{cases}
		\end{align*}
		the latter is valid, and so is \eqref{eq:BCQ}.
	\end{example}

	The final example of this section illustrates that, despite out above comments,
	\eqref{eq:BCQ} is not generally valid.
	\begin{example}\label{ex:BCQ_fails}
		Let us consider the parametric optimization problem
		\begin{equation}\label{eq:llp_BCQ_fails}
			\min\limits_y\{x(y_1+y_2)\,|\,y_1+y_2\leq 2,\,y_1-y_2\leq 0\}.
		\end{equation}
		For $x\in\R$ and $u\in\R^2_+$, we have
		\begin{align*}
			\psi_\ell(x,u)
			&=
			\inf_y\{x(y_1+y_2)+u_1(y_1+y_2-2)+u_2(y_1-y_2)\,|\,y\in\R^2\}
			\\
			&=
			\inf_y\{y_1(x+u_1+u_2)+y_2(x+u_1-u_2)-2u_1\,|\,y\in\R^2\}
			\\
			&=
			\begin{cases}
				-2u_1	&	x+u_1+u_2=0,\,x+u_1-u_2=0,
				\\
				-\infty	&	\text{otherwise.}
			\end{cases}
		\end{align*}
		Defining $\Omega\subset\R\times\R^2\times\R^2$ by means of
		\[
			\Omega
			:=
			\{(x,y,u)\in\R\times\R^2\times\R^2_+\,|\,x+u_1+u_2=0,\,x+u_1-u_2=0\},
		\]
		\cite[Exercise~10.10]{RockafellarWets1998} yields
		\[
			\partial^\infty\chi(x,y,u)=N_\Omega(x,y,u)
		\]
		for all $(x,y,u)\in\dom\chi=\Omega$ as we have
		\[
			\chi(x,y,u)
			=
			x(y_1+y_2)+2u_1+\delta_\Omega(x,y,u)
		\]
		in the present situation.
		Let us consider $(\bar x,\bar y,\bar u):=(0,(0,1),(0,0))$.
		Then a simple calculation reveals
		\begin{align*}
			N_\Omega(\bar x,\bar y,\bar u)
			&=
			\{(\xi_1+\xi_2,0,0,\xi_1+\xi_2+\xi_3,\xi_1-\xi_2+\xi_4)\in\R^5\,|\,\xi_3,\xi_4\leq 0\},
			\\
			N_{\mathcal F}(\bar x,\bar y,\bar u)
			&=
			\{(0,0,0)\}\times\R^2_-.
		\end{align*}
		Hence, for a vector $\eta\in\partial^\infty\chi(\bar x,\bar y,\bar u)\cap(-N_{\mathcal F}(\bar x,\bar y,\bar u))$,
		there exist $\xi_1,\xi_2\in\R$ and $\xi_3,\xi_4\in\R_-$ such that
		\[
			\eta_1=\xi_1+\xi_2=0,\,\eta_2=\eta_3=0,\,\eta_4=\xi_1+\xi_2+\xi_3\geq 0,\,\eta_5=\xi_1-\xi_2+\xi_4\geq 0.
		\]
		Choosing $\xi_1:=1$, $\xi_2:=-1$, $\xi_3:=\xi_4:=0$, we see that $\eta\neq 0$, 
		i.e., \eqref{eq:BCQ} is violated.\\
		It should be noted that despite the failure of \eqref{eq:BCQ},
		NSMFCQ fails to hold at $(\bar x,\bar y,\bar u)$ for the constraint system
		\[
			y_1+y_2\leq 2,\,y_1-y_2\leq 0,\,u\geq 0,\,x(y_1+y_2)-\psi_\ell(x,u)\leq 0,
		\]
		i.e., when considering some bilevel optimization problem with lower-level problem \eqref{eq:llp_BCQ_fails},
		then NSMFCQ fails for \eqref{eq:LDref} at $(\bar x,\bar y,\bar u)$
		(provided this point is feasible).		
	\end{example}

	It remains an open question whether NSMFCQ fails at all feasible points of \eqref{eq:LDref},
	i.e., even in the absence of \eqref{eq:BCQ}.

\subsection{Transformation based on Wolfe duality}\label{sec:trafo_Wolfe_duality}

Based on our considerations in \cref{sec:Wolfe_duality},
it seems reasonable to replace \eqref{eq:eoblp} by the single-level optimization problem
\begin{equation}\label{eq:WDref}\tag{WD$_\textup{ref}$}
  \begin{aligned}
  	&\min\limits_{x,y,z,u}&	&F(x,y)&		\\
  	&\subjectto&			&x\in X,&		\\
  	&&						&g(x,y)\leq 0,\,u\geq 0,&	\\
  	&&						&f(x,y)\leq L(x,z,u),&	\\
  	&&						&\nabla_2L(x,z,u)=0&
  \end{aligned}
\end{equation}
which we will refer to as the Wolfe dual reformulation of \eqref{eq:eoblp}.
This reformulation approach has been suggested in \cite{LiLinZhangZhu2022}.
Note that for each feasible point
$(\bar x,\bar y,\bar z,\bar u)\in\R^n\times\R^m\times\R^m\times\R^p$
of \eqref{eq:WDref}, \cref{prop:Wolfe_duality}\,\ref{item:Wolfe_weak}
ensures that $\bar y$ solves \hyperref[eq:llp]{\textup{(P$(\bar x)$)}}
to optimality while $(\bar z,\bar u)$ solves
\begin{equation}\label{eq:parametric_Wolfe_dual}
	\max\limits_{z,u}\{L(\bar x,z,u)\,|\,\nabla_2L(\bar x,z,u)=0,\,u\geq 0\}
\end{equation}
to global optimality.
Let us note that, in contrast to \eqref{eq:eoblp},
\eqref{eq:WDref} possesses the additional variables $(z,u)$.

\begin{remark}\label{rem:WDref_removal_z}
	Whenever the lower-level problem \eqref{eq:llp} is of the special form
	\[
		\min\limits_y\{f_1(x)+c^\top y\,|\,Ay\leq b(x)\}
	\]
	for continuously differentiable functions $f_1\colon\R^n\to\R$ and $b\colon\R^n\to\R^p$
	as well as some matrices $c\in\R^m$ and $A\in\R^{p\times m}$,
	\eqref{eq:WDref} takes the form
	\[	
		\begin{aligned}
  			&\min\limits_{x,y,z,u}&	&F(x,y)&		\\
  			&\subjectto&			&x\in X,&		\\
  			&&						&Ay\leq b(x),\,u\geq 0,&	\\
  			&&						&c^\top y+u^\top b(x)\leq 0,\,c+A^\top u=0&	
  		\end{aligned}
  	\]
  	and, thus, does not depend on the variable $z$.
  	Instead, the latter corresponds to the Lagrange dual reformulation \eqref{eq:LDref} 
  	of the original bilevel optimization problem as we have
  	\[
  		\psi_\ell(x,u)
  		=
  		\begin{cases}
  			f_1(x)-u^\top b(x)	& c+A^\top u=0,\\
  			-\infty				&\text{otherwise}
  		\end{cases}
  	\]
  	for each pair $(x,u)\in\R^n\times\R^p_+$ in the present situation.
\end{remark}

We are now concerned with the relationship between the global and local minimizers of
\eqref{eq:eoblp} and \eqref{eq:WDref}.
This topic already has been investigated in \cite[Section~3]{LiLinZhangZhu2022},
but we are going to address it from the viewpoint of the more general study \cite{BenkoMehlitz2021}.
Therefore, we make use of the intermediate mapping
$K_\textup{w}\colon\R^n\times\R^m\tto\R^m\times\R^p$ given by
\[
	\forall x\in\R^n,\,\forall y\in\R^m\colon\quad
	K_\textup{w}(x,y)
	:=
	\left\{
		(z,u)\in\R^m\times\R^p\,\middle|\,
		\begin{aligned}
			&g(x,y)\leq 0,\,u\geq 0,\\
			&f(x,y)\leq L(x,z,u),\\
			&\nabla_2L(x,z,u)=0
		\end{aligned}
	\right\}.
\]
Some properties of this mapping are listed in the following lemma.
\begin{lemma}\label{lem:properties_intermediate_map_Wolfe}
	\begin{enumerate}
		\item\label{item:K_w_closed} The set $\gph K_\textup{w}$ is closed.
		\item\label{item:K_w_dom1} Fix $(\bar x,\bar y)\in\gph\Psi$ and assume that there is some
			point in $\Psi(\bar x)$ where GCQ
			holds for \hyperref[eq:llp]{\textup{(P$(\bar x)$)}}.
			Then $(\bar x,\bar y)\in\dom K_\textup{w}$.
		\item\label{item:K_w_dom2} Fix $(\bar x,\bar y)\in\gph\Psi$.
			Then $(\bar y,\bar u)\in K_\textup{w}(\bar x,\bar y)$
			holds for each $\bar u\in\Lambda(\bar x,\bar y)$.
	\end{enumerate}
\end{lemma}
\begin{proof}
	The first assertion is obvious by continuous differentiability of
	the data functions $f$ and $g$.
	
	For the proof of the second assertion,
	suppose that GCQ holds at $z\in\Psi(\bar x)$
	for \hyperref[eq:llp]{\textup{(P$(\bar x)$)}}.
	Then, due to \cref{prop:Wolfe_duality}\,\ref{item:Wolfe_strong},
	there exists $u\in\R^p$ such that $(z,u)\in K_\textup{w}(\bar x,z)$.
	Due to $\bar y,z\in\Psi(\bar x)$, we have $f(\bar x,z)=f(\bar x,\bar y)$
	and $g(\bar x,\bar y)\leq 0$, yielding $(z,u)\in K_\textup{w}(\bar x,\bar y)$,
	i.e., $(\bar x,\bar y)\in\dom K_\textup{w}$.
	
	Let us now verify the third assertion.
	Therefore, we fix $\bar u\in\Lambda(\bar x,\bar y)$.
	This gives $\nabla_2 L(\bar x,\bar y,\bar u)=0$ and $\bar u\geq 0$.
	Furthermore, validity of the complementarity slackness condition also yields
	$f(\bar x,\bar y)=L(\bar x,\bar y,\bar u)$.
	Finally, $g(\bar x,\bar y)\leq 0$ is trivial from $\bar y\in\Psi(\bar x)$.
	Hence, $(\bar y,\bar u)\in K_\textup{w}(\bar x,\bar y)$ has been shown.
\end{proof}

The following result is similar to \cite[Theorem~4.3]{BenkoMehlitz2021}
and slightly generalizes \cite[Theorem~3.1]{LiLinZhangZhu2022} in the convex setting.
Note the similarities to the proof of \cref{thm:LDref_global}.

\begin{theorem}\label{thm:WDref_global}
	\begin{enumerate}
		\item\label{item:WDref_global}
			Let $(\bar x,\bar y)\in\R^n\times\R^m$ be a global minimizer of \eqref{eq:eoblp}.
			Then, for each $(\bar z,\bar u)\in K_\textup{w}(\bar x,\bar y)$,
			$(\bar x,\bar y,\bar z,\bar u)$ is a global minimizer of \eqref{eq:WDref}.
		\item Let $(\bar x,\bar y,\bar z,\bar u)\in\R^n\times\R^m\times\R^m\times\R^p$
			be a global minimizer of \eqref{eq:WDref}.
			Assume that, for each $x\in X\cap\dom\Psi$, there exists some point
			in $\Psi(x)$ where GCQ holds for \eqref{eq:llp}.
			Then $(\bar x,\bar y)$ is a global minimizer of \eqref{eq:eoblp}.
	\end{enumerate}
\end{theorem}
\begin{proof}
	We start with the proof of the first assertion.
	Let us assume that there is some $(\bar z,\bar u)\in K_\textup{w}(\bar x,\bar y)$
	such that $(\bar x,\bar y,\bar z,\bar u)$ is not a global minimizer of
	\eqref{eq:WDref}. Then there exists a feasible point
	$(\tilde x,\tilde y,\tilde z,\tilde u)\in\R^n\times\R^m\times\R^m\times\R^p$
	of \eqref{eq:WDref} such that $F(\bar x,\bar y)>F(\tilde x,\tilde y)$.
	As \cref{prop:Wolfe_duality}\,\ref{item:Wolfe_weak} yields $\tilde y\in\Psi(\tilde x)$,
	$(\tilde x,\tilde y)$ is a feasible point of \eqref{eq:eoblp},
	which means that $(\bar x,\bar y)$ cannot be a global minimizer of \eqref{eq:eoblp}.
	
	We now proceed with the proof of the second assertion.
	Due to \cref{prop:Wolfe_duality}\,\ref{item:Wolfe_weak},
	$(\bar x,\bar y)$ is feasible to \eqref{eq:eoblp}.
	If $(\bar x,\bar y)$ is not a global minimizer of \eqref{eq:eoblp},
	there exists some feasible point $(\tilde x,\tilde y)\in\R^n\times\R^m$ of \eqref{eq:eoblp}
	such that $F(\bar x,\bar y)>F(\tilde x,\tilde y)$.
	Due to $\tilde x\in X\cap\dom\Psi$ and the postulated assumption,
	\cref{lem:properties_intermediate_map_Wolfe}\,\ref{item:K_w_dom1} guarantees the existence of
	$(\tilde z,\tilde u)\in K_\textup{w}(\tilde x,\tilde y)$.
	Hence, $(\tilde x,\tilde y,\tilde z,\tilde u)$ is feasible to \eqref{eq:WDref}
	which shows that $(\bar x,\bar y,\bar z,\bar u)$ cannot be a global minimizer
	of that problem.
\end{proof}

It should be noted that, in the absence of a constraint qualification,
$K_\textup{w}(\bar x,\bar y)$ might be empty for some pair $(\bar x,\bar y)\in\gph\Psi$,
i.e., it is reasonable to assume validity of GCQ at
some point in $\Psi(\bar x)$ for \hyperref[eq:llp]{\textup{(P$(\bar x)$)}}
to make \cref{thm:WDref_global}\,\ref{item:WDref_global} meaningful,
see \cite[Theorem~3.1]{LiLinZhangZhu2022}.
For the overall transformation to be reasonable,
it is, thus, desirable that a constraint qualification holds at each feasible point
of the lower-level problem \eqref{eq:llp} for all $x\in X$.

Let us now investigate the relationship between local minimizers of
\eqref{eq:eoblp} and \eqref{eq:WDref}.
First, let us state a result which parallels \cite[Theorem~4.5]{BenkoMehlitz2021},
and its proof is similar to the one of \cref{thm:LDref_local}.
It is also loosely connected to \cite[Theorem~3.3]{LiLinZhangZhu2022} where
local optimality for \eqref{eq:WDref} is defined in such a way that neighborhoods
are only taken w.r.t.\ those variables which appear in the original bilevel
optimization problem \eqref{eq:eoblp}. This notion, however, is difficult to treat
from an algorithmic point of view and, thus, of limited practical use.

\begin{theorem}\label{thm:WDref_local}
	\begin{enumerate}
		\item\label{item:WDref_local_simple} Let $(\bar x,\bar y)\in\R^n\times\R^m$ be a local minimizer of \eqref{eq:eoblp}.
			Then, for each $(\bar z,\bar u)\in K_\textup{w}(\bar x,\bar y)$,
			$(\bar x,\bar y,\bar z,\bar u)$ is a local minimizer of \eqref{eq:WDref}.
		\item\label{item:WDref_local}
			Let $(\bar x,\bar y, z, u)\in\R^n\times\R^m\times\R^m\times\R^p$
			be a local minimizer of \eqref{eq:WDref} for each $(z,u)\in K_\textup{w}(\bar x,\bar y)$.
			Assume that, for each $x\in X\cap \dom\Psi$ in a neighborhood of $\bar x$,
			there exists some point in $\Psi(x)$ where GCQ holds for \eqref{eq:llp}.
			Finally, let $K_\textup{w}$ be inner semicompact at $(\bar x,\bar y)$ w.r.t.\ its domain.
			Then $(\bar x,\bar y)$ is a local minimizer of \eqref{eq:eoblp}.
	\end{enumerate}
\end{theorem}
\begin{proof}
	The proof of the first assertion follows in similar way as the first statement of \cref{thm:WDref_global}.
	We, thus, only prove the second statement.
	Due to \cref{prop:Wolfe_duality}\,\ref{item:Wolfe_weak},
	$(\bar x,\bar y)$ is feasible to \eqref{eq:eoblp}.
	Suppose that $(\bar x,\bar y)$ is not a local minimizer of \eqref{eq:eoblp}.
	Then we find a sequence $\{(x^k,y^k)\}_{k\in\N}\subset (X\times\R^m)\cap\gph\Psi$ of points
	such that $x^k\to\bar x$, $y^k\to\bar y$, and $F(x^k,y^k)<F(\bar x,\bar y)$ for all $k\in\N$.
	For sufficiently large $k\in\N$, \cref{lem:properties_intermediate_map_Wolfe}\,\ref{item:K_w_dom1}
	ensures $(x^k,y^k)\in\dom K_\textup{w}$.
	As $K_\textup{w}$ is assumed to be inner semicompact w.r.t.\ its domain at $(\bar x,\bar y)$,
	we can find a sequence $\{(z^k,u^k)\}_{k\in\N}\subset\R^m\times\R^p$ such that
	$(z^k,u^k)\in K_\textup{w}(x^k,y^k)$ holds for all large enough $k\in\N$
	and $\{(z^k,u^k)\}_{k\in\N}$ possesses a bounded subsequence.
	Along a subsequence (without relabeling), we may assume $z^k\to\bar z$ and $u^k\to\bar u$
	for a pair $(\bar z,\bar u)\in\R^m\times\R^p$.
	Applying \cref{lem:properties_intermediate_map_Wolfe}\,\ref{item:K_w_closed},
	$(\bar z,\bar u)\in K_\textup{w}(\bar x,\bar y)$ follows, i.e.,
	$(\bar x,\bar y,\bar z,\bar u)$ is feasible to \eqref{eq:WDref}.
	Since all points along the tail of the sequence $\{(x^k,y^k,z^k,u^k)\}_{k\in\N}$ possess a smaller
	objective value than $(\bar x,\bar y,\bar z,\bar u)$,
	the latter cannot be a local minimizer of \eqref{eq:WDref}.	
	This, however, contradicts the theorem's assumptions.
\end{proof}

Subsequently, we reinspect the classical example from \cite[Example~3.4]{DempeDutta2012}
in order to illustrate that, in general, local minimizers of \eqref{eq:WDref} do not
correspond to local minimizers of \eqref{eq:eoblp}.
Thereby, we extend the findings from \cref{ex:artificial_local_solution_l} to the
single-level reformulation \eqref{eq:WDref}.

\begin{example}\label{ex:artificial_local_solution_w}
	We consider the bilevel optimization problem \eqref{eq:ulp_ex_lagrange}
	with associated lower-level problem \eqref{eq:llp_ex_lagrange}
	which we already discussed in \cref{ex:artificial_local_solution_l}.

	The associated problem \eqref{eq:WDref} reads as
	\[
		\begin{aligned}
		&\min\limits_{x,y,z,u}&	&(x-1)^2+(y-1)^2&	\\
		&\subjectto&			& x+y\leq 1,\,-x+y\leq 1,\, u\geq 0,&\\
		&&						& -y \leq -z + u_1(x+z-1) + u_2(-x+z-1),&\\
		&&						& -1 + u_1 + u_2 = 0,&
		\end{aligned}
	\]
	and the latter can be simplified to
	\begin{equation}\label{eq:WDref_ex_Wolfe1}
		\begin{aligned}
		&\min\limits_{x,y,z,u}&	&(x-1)^2+(y-1)^2&	\\
		&\subjectto&			& x+y\leq 1,\,-x+y\leq 1,\, u\geq 0,&\\
		&&						& y \geq 1-(u_1-u_2)x,&\\
		&&						& u_1 + u_2 = 1,&
		\end{aligned}
	\end{equation}
	see \cref{rem:WDref_removal_z} as well.
	We note that, apart from the fact that, in \eqref{eq:WDref_ex_Wolfe1},
	one also has to minimize over the variable $z$ which has been removed from
	the problem, \eqref{eq:WDref_ex_Wolfe1} is the same problem as
	\eqref{eq:LDref_example}.
	Hence, due to our considerations in \cref{ex:artificial_local_solution_l},
	the point $(\bar x,\bar y,\bar z,\bar u):=(0,1,1,(0,1))$,
	which is feasible to \eqref{eq:WDref_ex_Wolfe1},
	is a local minimizer of that problem
	while $(\bar x,\bar y)$ is not a local minimizer of \eqref{eq:ulp_ex_lagrange}.
	
	We emphasize that $K_\textup{w}(\bar x,\bar y)=\{(z,u)\in\R\times\R^2_+\,|\,u_1+u_2=1\}$
	holds true.
	Hence, for $\tilde u:=(1,0)$, we have $(\bar z,\tilde u)\in K_\textup{w}(\bar x,\bar y)$.
	Furthermore, for each $\varepsilon>0$,
	$(\bar x+\varepsilon,\bar y-\varepsilon,\bar z,\tilde u)$
	is feasible to \eqref{eq:WDref_ex_Wolfe1} and, for $\varepsilon\in(0,1)$, possesses a better
	objective value than $(\bar x,\bar y,\bar z,\tilde u)$.
	Hence, $(\bar x,\bar y,\bar z,\tilde u)$ is not a local minimizer
	of \eqref{eq:WDref_ex_Wolfe1}.
	That is why it is not possible to abstain from demanding
	local optimality of $(\bar x,\bar y,z,u)$
	for all $(z,u)\in K_\textup{w}(\bar x,\bar y)$
	in \cref{thm:WDref_local}\,\ref{item:WDref_local}.
	Let us hint the interested reader to the fact that all the other assumptions in
	\cref{thm:WDref_local}\,\ref{item:WDref_local} are satisfied.
	Validity of GCQ for all feasible points of the lower-level problem \eqref{eq:llp_ex_lagrange}
	has been discussed in \cref{ex:artificial_local_solution_l} already.
	Furthermore, one can easily check that $K_\textup{w}$ is inner semicompact at each point
	of its domain w.r.t.\ its domain since the variable $z$ is free while
	$u$ has to be chosen from the (bounded) standard simplex $\{u\in\R^2_+\,|\,u_1+u_2=1\}$
	(subject to yet another linear constraint).
	We also want to emphasize that, due to the fact that the variable $z$ can be chosen
	freely in the images of $K_\textup{w}$, this mapping possesses unbounded images.
	This, however, does not effect the necessary inner semicompactness.
\end{example}

Let us now take a closer look at the rather abstract inner semicompactness requirement
on the intermediate mapping $K_\textup{w}$ which has been exploited in
\cref{thm:WDref_local}\,\ref{item:WDref_local}.
As mentioned earlier,
for some feasible point $(\bar x,\bar y,\bar z,\bar u)\in\R^n\times\R^m\times\R^m\times\R^p$
of \eqref{eq:WDref}, we know $\bar y\in\Psi(\bar x)$ and that $(\bar z,\bar u)$ is a
global maximizer of \eqref{eq:parametric_Wolfe_dual}. Furthermore, due to $\bar u\geq 0$ and
our inherent convexity assumptions, we know that $\nabla^2_{2,2}L(\bar x,\bar z,\bar u)$ is
positive semidefinite. Let us assume for a moment that $\nabla^2_{2,2}L(\bar x,\cdot,\bar u)$ is
uniformly positive definite on $\R^m$.
This can be ensured if $f(\bar x,\cdot)$ is uniformly convex or if there is some $i\in\{1,\ldots,p\}$ such that
$\bar u_i>0$ and $g_i(\bar x,\cdot)$ is uniformly convex.
We also note that these assumptions are stable under small perturbations of $\bar x$ (in $\R^n$) and $\bar u$ (in $\R^p_+$).
Furthermore, in \cref{thm:WDref_local}\,\ref{item:WDref_local}, we anyhow require that,
for all $x\in X\cap\dom\Psi$ locally around $\bar x$, there is a point in $\Psi(x)$ where
GCQ holds. Hence, we can apply \cref{prop:Wolfe_duality}\,\ref{item:Wolfe_strong_converse}
in order to find that, for each feasible point $(x,y,z,u)\in\R^n\times\R^m\times\R^m\times\R^p$
of \eqref{eq:WDref} close to $(\bar x,\bar y,\bar z,\bar u)$, we have $z=y$ and $u\in\Lambda(x,y)$.
Hence, the inner semicompactness requirement in \cref{thm:WDref_local}\,\ref{item:WDref_local}
reduces to the one of $\Lambda$, and the latter already popped up in \cref{sec:trafo_Lagrange_duality},
see the discussion following the proof of \cref{thm:LDref_local}.

This observation motivates the subsequently stated theorem.
We emphasize that, in contrast to the arguments drafted above,
we do not exploit any additional properties of the data functions $f$ and $g$,
as we do not rely on \cref{prop:Wolfe_duality}\,\ref{item:Wolfe_strong_converse} for its proof.
However, a regularity condition stronger than GCQ is needed.
This observation can also be found in \cite[Theorem~3.4]{LiLinZhangZhu2022}.

\begin{theorem}\label{thm:WDref_local_refined}
	Let $(\bar x,\bar y, \bar y, u)\in\R^n\times\R^m\times\R^m\times\R^p$
	be a local minimizer of \eqref{eq:WDref} for each $u\in\Lambda(\bar x,\bar y)$.
	Furthermore, assume that Slater's constraint qualification
	holds for \hyperref[eq:llp]{\textup{(P$(\bar x)$)}}.
	Then $(\bar x,\bar y)$ is a local minimizer of \eqref{eq:eoblp}.
\end{theorem}
\begin{proof}
	Due to \cref{prop:Wolfe_duality}\,\ref{item:Wolfe_weak},
	$(\bar x,\bar y)$ is feasible to \eqref{eq:eoblp}.
	Let us now proceed as in the proof of \cref{thm:WDref_local}\,\ref{item:WDref_local}.
	We suppose that $(\bar x,\bar y)$ is not a local minimizer of \eqref{eq:eoblp}.
	Then we find a sequence $\{(x^k,y^k)\}_{k\in\N}\subset(X\times\R^m)\cap\gph\Psi$
	such that $x^k\to\bar x$, $y^k\to\bar y$, and $F(x^k,y^k)<F(\bar x,\bar y)$ for all $k\in\N$.
	Validity of Slater's constraint qualification for \hyperref[eq:llp]{\textup{(P$(\bar x)$)}}
	implies validity of Slater's constraint qualification
	for \hyperref[eq:llp]{\textup{(P$(x^k)$)}} provided $k\in\N$ is sufficiently large.
	This particularly guarantees that MFCQ and, thus,
	GCQ hold at each feasible point of these problems.
	Hence, we can pick $u^k\in\Lambda(x^k,y^k)$ for sufficiently large $k\in\N$,
	and \cref{lem:properties_intermediate_map_Wolfe}\,\ref{item:K_w_dom2} guarantees
	$(y^k,u^k)\in K_\textup{w}(x^k,y^k)$.
	Consulting \cite[Proposition~4.43]{BonnansShapiro2000} or \cite[proof of Theorem~4.2]{Dempe2002}, 
	the Lagrange multiplier mapping
	$\Lambda$ possesses uniformly bounded image sets in a neighborhood of $(\bar x,\bar y)$.
	Consequently, along a subsequence (without relabeling), there is some $\bar u\in\R^p$
	such that $u^k\to\bar u$,
	and the closedness of $\gph\Lambda$ yields $\bar u\in\Lambda(\bar x,\bar y)$.
	By \cref{lem:properties_intermediate_map_Wolfe}\,\ref{item:K_w_dom2},
	the point $(x^k,y^k,y^k,u^k)$ is feasible to \eqref{eq:WDref} for each $k\in\N$
	and possesses a smaller objective value than $(\bar x,\bar y,\bar y,\bar u)$ which is
	also feasible to \eqref{eq:WDref}, i.e., the latter cannot be a local minimizer of
	\eqref{eq:WDref}, contradicting the theorem's assumptions.
\end{proof}

Note the similarities between the proofs of \cref{thm:WDref_local}\,\ref{item:WDref_local}
and \cref{thm:WDref_local_refined}. In fact, in the proof of \cref{thm:WDref_local_refined},
we show that the abstract inner semicompactness of $K_\textup{w}$ which is assumed in
\cref{thm:WDref_local}\,\ref{item:WDref_local} is (exemplary) realized by pairs
$(y^k,u^k)\in K_\textup{w}(x^k,y^k)$, $k\in\N$, such that $u^k\in\Lambda(x^k,y^k)$
thanks to \cref{lem:properties_intermediate_map_Wolfe}\,\ref{item:K_w_dom2}.
This procedure also allows to restrict to feasible points of \eqref{eq:WDref}
whose $y$- and $z$-component are the same, which makes \cref{thm:WDref_local_refined} much
more tractable than \cref{thm:WDref_local}\,\ref{item:WDref_local}.
However, one is in need of stronger qualification conditions in \cref{thm:WDref_local_refined}.
Again, \cref{ex:artificial_local_solution_w} illustrates that we cannot get rid of the
assumption that local optimality of $(\bar x,\bar y,\bar y,u)$ for \eqref{eq:WDref}
has to holds for all $u\in\Lambda(\bar x,\bar y)$.
Observe that \cref{thm:WDref_local_refined} is a perfect companion of
\cite[Theorem~3.2]{DempeDutta2012} which addresses \eqref{eq:KKTref}.

As mentioned in \cref{sec:trafo_Lagrange_duality}, one can rely on qualification conditions
different from Slater's constraint qualification in order to guarantee inner semicompactness
of the multiplier mapping $\Lambda$ w.r.t.\ its domain in \cref{thm:WDref_local_refined}.
Then, however, one potentially has to ensure regularity of the problem \eqref{eq:llp} for all $x\in X\cap\dom\Psi$
locally around $\bar x$ at its minimizers via additional assumptions.

It has been shown in \cite[Example~4.1]{LiLinZhangZhu2022}
by means of an example that MFCQ is not necessarily
violated at the feasible points of the surrogate optimization problem
\eqref{eq:WDref}. This example, however, is awkwardly misleading as the chosen
lower-level problem in \cite[Example~4.1]{LiLinZhangZhu2022} does not
satisfy weak Wolfe duality, which also means that the reformulation
\eqref{eq:WDref} is not at all reasonable.
Indeed, the authors in \cite{LiLinZhangZhu2022} investigate the
lower-level problem
\[
	\min\limits_y\{y\,|\,y^3\leq x,\,y\geq 0\}.
\]
The associated Wolfe dual problem reads as
\[
	\max\limits_{z,u}\{z+u_1(z^3-x)-u_2z\,|\,u\geq 0,\,1+3u_1z^2-u_2=0\},
\]
and inserting $u_2=1+3u_1z^2$ into the objective function, the latter can be stated equivalently in the form
\[
	\max\limits_{z,u}\{-2u_1z^3-u_1x\,|\,u\geq 0,\,1+3u_1z^2-u_2=0\}.
\]
Note that, for each $x\geq 0$, $\bar y:=0$ is feasible for the primal problem and possesses objective value $0$,
while, e.g., $(\bar z,\bar u_\varepsilon):=(-3,(\varepsilon,1+27\varepsilon))$ for arbitrary $\varepsilon>0$,
is feasible for the Wolfe dual problem and possesses objective value $(54-x)\varepsilon$.
The latter is strictly positive for each $x<54$,
showing that weak Wolfe duality fails for all $x\in[0,54)$.
In \cite[Example~4.1]{LiLinZhangZhu2022}, the authors considered the point
$(\bar x,\bar y,\bar z,\bar u):=(8,0,-3,(0,1))$ so that the above is a relevant issue.

In the following theorem, we indicate that, in the convex setting we are considering here,
MFCQ is violated at each feasible point of \eqref{eq:WDref}.
The proof is similar to the one of \cref{thm:violation_MFCQ_Lagrange}.
When transferring this result to more general situations, it still remains true as long as
the lower-level problem \eqref{eq:llp} satisfies weak Wolfe duality for each $x\in X$,
which can be easily observed from its proof.
\begin{theorem}\label{thm:violation_of_MFCQ_Wolfe}
	Let $(\bar x,\bar y,\bar z,\bar u)\in\R^n\times\R^m\times\R^m\times\R^p$ be a feasible point of \eqref{eq:WDref}.
	Furthermore, assume that $X:=\{x\in\R^n\,|\,G(x)\leq 0\}$ holds
 	for a continuously differentiable function $G\colon\R^n\to\R^q$.
 	Then MFCQ for \eqref{eq:WDref} is violated at $(\bar x,\bar y,\bar z,\bar u)$.
\end{theorem}
\begin{proof}
	If MFCQ fails to hold for the constraint system
	\[
		g(x,y)\leq 0,\,u\geq 0,\,\nabla_2L(x,z,u)=0
	\]
	at $(\bar x,\bar y,\bar z,\bar u)$, then it is also violated for \eqref{eq:WDref} at the same point.
	Thus, let us assume that MFCQ holds for the above constraint system at $(\bar x,\bar y,\bar z,\bar u)$.
	Note that \cref{prop:Wolfe_duality}\,\ref{item:Wolfe_weak} yields
	$f(\bar x,\bar y)-L(\bar x,\bar z,\bar u)=0$, and that $(\bar x,\bar y,\bar z,\bar u)$ is a minimizer
	of the optimization problem
	\begin{equation}\label{eq:violation_MFCQ_Wolfe_surrogate}
		\min\limits_{x,y,z,u}\{f(x,y)-L(x,z,u)\,|\,g(x,y)\leq 0,\,u\geq 0,\,\nabla_2L(x,z,u)=0\}.
	\end{equation}
	Since MFCQ holds for this problem at $(\bar x,\bar y,\bar z,\bar u)$,
	the latter point is stationary, i.e., we find multipliers
	$\alpha,\beta\in\R^p$ and $\delta\in\R^m$ such that
	\begin{align*}
		\nabla_1f(x,y)-\nabla_1L(x,z,u)+g'_1(x,y)^\top\alpha+\nabla^2_{2,1}L(x,z,u)\delta&=0,
		\\
		\nabla_2f(x,y)+g'_2(x,y)^\top\alpha&=0,
		\\
		-\nabla_2L(x,z,u)+\nabla^2_{2,2}L(x,z,u)\delta&=0,
		\\
		-g(x,z)-\beta +g'_2(x,z)\delta&=0,
		\\
		\alpha\geq 0,\,\alpha^\top g(x,y)&=0,
		\\
		\beta\geq 0,\,\beta^\top u&=0.
	\end{align*}
	This, however, means that MFCQ fails
	for \eqref{eq:WDref} at $(\bar x,\bar y,\bar z,\bar u)$.
\end{proof}

The authors in \cite{LiLinZhangZhu2022} emphasize that a positive feature of the reformulation
\eqref{eq:WDref} is the fact that MFCQ may hold at some feasible points as this allows for
the application of standard solvers from nonlinear programming to solve \eqref{eq:WDref}.
Taking \cref{thm:violation_of_MFCQ_Wolfe} and our antedated investigation of
\cite[Example~4.1]{LiLinZhangZhu2022} as well as the associated discussion into account,
the only situation where MFCQ may hold at feasible points of \eqref{eq:WDref} is the one
where the lower-level problem does not satisfy weak Wolfe duality,
in which case the overall approach is nonsense.
For the purpose of completeness, let us hint the reader where the problem considered
in \cite[Example~4.1]{LiLinZhangZhu2022} eludes the arguments in the proof of
\cref{thm:violation_of_MFCQ_Wolfe}.
Due to the lack of weak Wolfe duality, the considered point $(\bar x,\bar y,\bar z,\bar u):=(8,0,-3,(0,1))$
is not a (local) minimizer of the associated problem \eqref{eq:violation_MFCQ_Wolfe_surrogate}.
Indeed, one can check that, for each $\varepsilon>0$, the point
$(\bar x,\bar y,\bar z,\bar u_\varepsilon)$ with $\bar u_\varepsilon:=(\varepsilon,1+27\varepsilon)$
is feasible to \eqref{eq:violation_MFCQ_Wolfe_surrogate} and possesses objective value
$-46\varepsilon<0$.
Thus, we cannot rely on the KKT conditions of \eqref{eq:violation_MFCQ_Wolfe_surrogate}
to verify the violation of MFCQ for \eqref{eq:WDref}.

In \cite{DiehlHouskaSteinSteuermann2013}, 
the authors suggest to reformulate a generalized semi-infinite optimization problem
based on Wolfe duality,
and it is shown that a problem-tailored version of MFCQ can hold at 
the feasible points of this reformulated problem.
Starting from \eqref{eq:VFR}, \eqref{eq:eoblp} obviously coincides with
\begin{equation}\label{eq:eoblp_as_gsip}
	\min\limits_{x,y}\{F(x,y)\,|\,(x,y)\in \mathcal S\cap\mathcal T\},
\end{equation}
where the sets $\mathcal S,\mathcal T\subset\R^n\times\R^m$ are given by
\begin{align*}
	\mathcal S
	:=&
	\{(x,y)\in X\times\R^m\,|\,g(x,y)\leq 0\},
	\\
	\mathcal T
	:=&
	\{(x,y)\in \R^n\times\R^m\,|\,f(x,y)\leq\varphi(x)\}
	\\
	=&
	\{(x,y)\in\R^n\times\R^m\,|\,\forall z\in \Gamma(x)\colon\,f(x,y)\leq f(x,z)\}.
\end{align*}
Above, we made use of the lower-level feasibility mapping $\Gamma$ 
from \eqref{eq:feasibility_mapping}.
The latter representation of $\mathcal T$ then induces the structure of
a generalized semi-infinite problem in \eqref{eq:eoblp_as_gsip},
see e.g.\ \cite[Section~3.1]{ZemkohoZhou2021} where this has been observed as well.
Following the arguments in \cite{DiehlHouskaSteinSteuermann2013}, 
it is reasonable to replace \eqref{eq:eoblp_as_gsip} by
\begin{equation}\label{eq:WD_ref_Stein}
	\begin{aligned}
		&\min\limits_{x,y,z,u}&		&F(x,y)&
		\\
		&\subjectto&				&(x,y)\in\mathcal S,\,u\geq 0,&
		\\
		&&							&f(x,y)\leq L(x,z,u),&
		\\
		&&							&\nabla_2L(x,z,u)=0,&
	\end{aligned}
\end{equation}
i.e., \eqref{eq:WDref}.
Taking a closer look at \cite{DiehlHouskaSteinSteuermann2013},
the authors show that \eqref{eq:WD_ref_Stein} may satisfy MFCQ
if $\mathcal S$ equals the whole space $\R^n\times\R^m$.
However, as $\mathcal S$ comprises the lower-level constraints,
this choice is impossible in the particular setting of
bilevel optimization as these constraints are essential in order
to connect \eqref{eq:WD_ref_Stein} with \eqref{eq:eoblp}.
Hence, \cref{thm:violation_of_MFCQ_Wolfe} does not contradict
the results in \cite{DiehlHouskaSteinSteuermann2013},
the latter are simply not applicable here.

Summing up the investigations in this subsection,
considering \eqref{eq:WDref} comes along with the same difficulties
that have to faced when \eqref{eq:KKTref} and \eqref{eq:LDref} are exploited.
To circumvent the potential lack of regularity,
one can rely on relaxation approaches, see e.g.\ \cite{DempeFranke2019}.
However, this does not resolve the issues related to the presence
of the artificial variables $(z,u)$ which generate artificial
stationary points, see e.g.\ \cref{ex:artificial_local_solution_w}.

\subsection{Transformation based on Mond--Weir duality}\label{sec:trafo_Mond_Weir_duality}

Invoking the theory from \cref{sec:Mond_Weir_duality},
another reasonable single-level reformulation of \eqref{eq:eoblp} is
\begin{equation}\label{eq:MWDref}\tag{MWD$_\textup{ref}$}
  \begin{aligned}
  	&\min\limits_{x,y,z,u}&	&F(x,y)&		\\
  	&\subjectto&			&x\in X,&		\\
  	&&						&g(x,y)\leq 0,\,u\geq 0,&	\\
  	&&						&f(x,y)\leq f(x,z),\,u^\top g(x,z)\geq 0,&	\\
  	&&						&\nabla_2L(x,z,u)=0&
  \end{aligned}
\end{equation}
which will be referred to as the Mond--Weir dual reformulation of \eqref{eq:eoblp}.
This transformation approach has been investigated in \cite{LiLinZhu2024} first.
Recalling \cref{prop:MondWeir_duality}\,\ref{item:MondWeir_weak},
for a feasible point $(\bar x,\bar y,\bar z,\bar u)\in\R^n\times\R^m\times\R^n\times\R^p$
of \eqref{eq:MWDref}, we know that $\bar y$ solves \hyperref[eq:llp]{\textup{(P$(\bar x)$)}}
to optimality while $(\bar z,\bar u)$ solves
\[
	\max\limits_{z,u}\{f(x,z)\,|\,\nabla_2L(x,z,u)=0,\,u^\top g(x,z)\geq 0,\,u\geq 0\}
\]
to global optimality.

Similar as in \cref{sec:trafo_Wolfe_duality}, we are going to analyze the relationship
between the global and local minimizers of \eqref{eq:eoblp} and \eqref{eq:MWDref},
deepening the analysis which can be found in \cite[Section~3]{LiLinZhu2024}.
Again, we follow \cite{BenkoMehlitz2021}, interpreting the variables $(z,u)$ in
\eqref{eq:MWDref} as implicit ones which have been made explicit.
Therefore, we introduce the intermediate mapping
$K_\textup{mw}\colon\R^n\times\R^m\tto\R^m\times\R^p$ by means of
\[
	\forall x\in\R^n,\,\forall y\in\R^m\colon\quad
	K_\textup{mw}(x,y)
	:=
	\left\{
		(z,u)\in\R^m\times\R^p\,\middle|\,
		\begin{aligned}
			&g(x,y)\leq 0,\,u\geq 0,\\
			&f(x,y)\leq f(x,z),\\
			&u^\top g(x,z)\geq 0,\\
			&\nabla_2L(x,z,u)=0
		\end{aligned}
	\right\}.
\]
Similarly as in \cref{lem:properties_intermediate_map_Wolfe},
we collect some properties of the intermediate mapping $K_\textup{mw}$
in the subsequently stated lemma.
The proof is analogous to the one of \cref{lem:properties_intermediate_map_Wolfe} and,
thus, omitted.

\begin{lemma}\label{lem:properties_intermediate_map_Mond_Weir}
	\begin{enumerate}
		\item The set $\gph K_\textup{mw}$ is closed.
		\item Fix $(\bar x,\bar y)\in\gph\Psi$ and assume that there is some point in
			$\Psi(\bar x)$ where GCQ holds.
			Then $(\bar x,\bar y)\in\dom K_\textup{mw}$.
		\item Fix $(\bar x,\bar y)\in\gph\Psi$.
			Then $(\bar y,\bar u)\in K_\textup{mw}(\bar x,\bar y)$ holds for each
			$\bar u\in\Lambda(\bar x,\bar y)$.
	\end{enumerate}
\end{lemma}

Let us start our investigations with some first observations regarding the relationship
between the single-level reformulations \eqref{eq:WDref} and \eqref{eq:MWDref}.
For each pair $(x,y)\in\R^n\times\R^m$, the inclusion
$K_\textup{mw}(x,y)\subset K_\textup{w}(x,y)$ is valid.
Particularly, we have $\gph K_\textup{mw}\subset\gph K_\textup{w}$.
Since the feasible set of \eqref{eq:WDref} is given by
\[
	(X\times\R^m\times\R^m\times\R^p)\cap\gph K_\textup{w}
\]
while the feasible set of \eqref{eq:MWDref} can be represented as
\[
	(X\times\R^m\times\R^m\times\R^p)\cap\gph K_\textup{mw},
\]
we obtain that each global (local) minimizer of \eqref{eq:WDref},
which is feasible to \eqref{eq:MWDref}, is already a global (local)
minimizer of \eqref{eq:MWDref}.
Let us also note that, in \cite[Theorem~5.1]{LiLinZhu2024},
it has been shown that a KKT point of \eqref{eq:WDref},
which is feasible to \eqref{eq:MWDref}, is already a KKT point of \eqref{eq:MWDref}.
Hence, \eqref{eq:MWDref} provides a reformulation of \eqref{eq:eoblp} which is, in general,
tighter than \eqref{eq:WDref}, at least w.r.t.\ the involved surrogate variables.

We are well-prepared for a comparison of the global and local minimizers associated
with \eqref{eq:eoblp} and \eqref{eq:MWDref}.
The following three theorems can,
with the aid of \cref{lem:properties_intermediate_map_Mond_Weir},
be proven as \cref{thm:WDref_global,thm:WDref_local,thm:WDref_local_refined},
so we abstain from a validation.
Related statement can be found in \cite[Theorem~3.1, 3.3, and~3.4]{LiLinZhu2024}.
We note, however, that \cref{thm:MWDref_local} significantly differs from
\cite[Theorem~3.3]{LiLinZhu2024} where, again, local minimality in \eqref{eq:MWDref}
is defined in a strange way which is difficult to exploit algorithmically.
Furthermore, \cref{thm:MWDref_local_refined} can be motivated via
\cref{prop:MondWeir_duality}\,\ref{item:MondWeir_strong_converse} in similar way
as \cref{thm:WDref_local_refined} has been motivated via
\cref{prop:Wolfe_duality}\,\ref{item:Wolfe_strong_converse}.

\begin{theorem}\label{thm:MWDref_global}
	\begin{enumerate}
		\item
			Let $(\bar x,\bar y)\in\R^n\times\R^m$ be a global minimizer of \eqref{eq:eoblp}.
			Then, for each $(\bar z,\bar u)\in K_\textup{mw}(\bar x,\bar y)$,
			$(\bar x,\bar y,\bar z,\bar u)$ is a global minimizer of \eqref{eq:MWDref}.
		\item Let $(\bar x,\bar y,\bar z,\bar u)\in\R^n\times\R^m\times\R^m\times\R^p$
			be a global minimizer of \eqref{eq:MWDref}.
			Assume that, for each $x\in X\cap\dom\Psi$, there exists some point
			in $\Psi(x)$ where GCQ holds for \eqref{eq:llp}.
			Then $(\bar x,\bar y)$ is a global minimizer of \eqref{eq:eoblp}.
	\end{enumerate}
\end{theorem}

\begin{theorem}\label{thm:MWDref_local}
	\begin{enumerate}
		\item Let $(\bar x,\bar y)\in\R^n\times\R^m$ be a local minimizer of \eqref{eq:eoblp}.
			Then, for each $(\bar z,\bar u)\in K_\textup{mw}(\bar x,\bar y)$,
			$(\bar x,\bar y,\bar z,\bar u)$ is a local minimizer of \eqref{eq:MWDref}.
		\item\label{item:MWDref_local}
			Let $(\bar x,\bar y, z, u)\in\R^n\times\R^m\times\R^m\times\R^p$
			be a local minimizer of \eqref{eq:MWDref}
			for each $(z,u)\in K_\textup{mw}(\bar x,\bar y)$.
			Assume that, for each $x\in X\cap \dom\Psi$ in a neighborhood of $\bar x$,
			there exists some point in $\Psi(x)$ where
			GCQ holds for \eqref{eq:llp}.
			Finally, let $K_\textup{mw}$ be inner semicompact at $(\bar x,\bar y)$
			w.r.t.\ its domain.
			Then $(\bar x,\bar y)$ is a local minimizer of \eqref{eq:eoblp}.
	\end{enumerate}
\end{theorem}

\begin{theorem}\label{thm:MWDref_local_refined}
	Let $(\bar x,\bar y, \bar y, u)\in\R^n\times\R^m\times\R^m\times\R^p$
	be a local minimizer of \eqref{eq:MWDref} for each $u\in\Lambda(\bar x,\bar y)$.
	Furthermore, assume that Slater's constraint qualification
	holds for \hyperref[eq:llp]{\textup{(P$(\bar x)$)}}.
	Then $(\bar x,\bar y)$ is a local minimizer of \eqref{eq:eoblp}.
\end{theorem}

The following example, which extends \cref{ex:artificial_local_solution_w},
illustrates that local minimizers of \eqref{eq:MWDref} do not correspond to
local minimizers of \eqref{eq:eoblp} in general.

\begin{example}\label{ex:artifical_local_solutions_mw}
	We reconsider the bilevel optimization problem \eqref{eq:ulp_ex_lagrange}
	with lower-level problem \eqref{eq:llp_ex_lagrange} which already has been
	investigated in \cref{ex:artificial_local_solution_l,ex:artificial_local_solution_w}.
	
	The associated problem \eqref{eq:MWDref} is given by
	\begin{equation}\label{eq:MWDref_ex}
		\begin{aligned}
		&\min\limits_{x,y,z,u}&	&(x-1)^2+(y-1)^2&	\\
		&\subjectto&			& x+y\leq 1,\,-x+y\leq 1,\, u\geq 0,&\\
		&&						& -y \leq -z,\,u_1(x+z-1)+u_2(-x+z-1)\geq 0,&\\
		&&						& -1 + u_1 + u_2 = 0.&
		\end{aligned}
	\end{equation}
	Let us note that $(\bar x,\bar y,\bar z,\bar u):=(0,1,1,(0,1))$ is feasible to \eqref{eq:MWDref_ex}.
	We emphasize that each feasible point $(x,y,z,u)\in\R\times\R\times\R\times\R^2$ of \eqref{eq:MWDref_ex}
	necessarily satisfies $z\leq y\leq 1$. This immediately yields
	\begin{equation}\label{eq:some_estimate}
		\begin{aligned}
		u_1(x+y-1)+u_2(-x+y-1)
		&
		=
		x(u_1-u_2)+y(u_1+u_2)-(u_1+u_2)
		\\
		&=
		x(u_1-u_2)+y-1
		\\
		&\geq
		x(u_1-u_2)+z-1
		\\
		&=
		x(u_1-u_2)+z(u_1+u_2)-(u_1+u_2)
		\\
		&=
		u_1(x+z-1)+u_2(-x+z-1)
		\geq
		0,
		\end{aligned}
	\end{equation}
	and due to the remaining constraints, this yields
	\[
		u_1(x+y-1)=0,\,u_2(-x+y-1)=0.
	\]
	Whenever $(x,y,z,u)$ is, additionally, chosen from a sufficiently small neighborhood of $(\bar x,\bar y,\bar z,\bar u)$, 
	we have $u_1\geq 0$ and $u_2>0$.
	If $u_1>0$, the above yields $x=0=\bar x$ and $y=1=\bar y$, i.e., the objective value of $(x,y,z,u)$
	is the same as the one for $(\bar x,\bar y,\bar z,\bar u)$.
	In the other case $u_1=0$, we have $u_1-u_2<0$. Due to $z\leq 1$, we find $x(u_1-u_2)\geq 1-z\geq 0$
	from \eqref{eq:some_estimate}, and this yields $x\leq 0$.
	Then, however, the objective value of $(x,y,z,u)$ cannot be smaller than the one of $(\bar x,\bar y,\bar z,\bar u)$.
	Hence, the latter is a local minimizer of \eqref{eq:MWDref_ex}
	while $(\bar x,\bar y)$ is not a local minimizer of \eqref{eq:ulp_ex_lagrange}.
	
	A simple calculation reveals
	$K_\textup{mw}(\bar x,\bar y)=\{(1,u)\in\R\times\R^2_+\,|\,u_1+u_2=1\}=\{1\}\times\Lambda(\bar x,\bar y)$.
	As in \cref{ex:artificial_local_solution_l,ex:artificial_local_solution_w},
	one can easily check that, for $\tilde u:=(1,0)$,
	$(\bar x,\bar y,\bar z,\tilde u)$ is not a local minimizer of \eqref{eq:MWDref_ex}.
	This means that one cannot abstain from postulating local optimality
	of $(\bar x,\bar y,z,u)$ for all $(z,u)\in K_\textup{mw}(\bar x,\bar y)$
	in \cref{thm:MWDref_local}\,\ref{item:MWDref_local} and local optimality
	of $(\bar x,\bar y,\bar y,u)$ for all $u\in\Lambda(\bar x,\bar y)$
	in \cref{thm:MWDref_local_refined},
	as all the other assumptions of \cref{thm:MWDref_local}\,\ref{item:MWDref_local}
	and \cref{thm:MWDref_local_refined} are clearly valid.
	Note that the inner semicompactness of the intermediate mapping $K_\textup{mw}$ w.r.t.\ its domain
	at $(\bar x,\bar y)$ follows from the representation
	\[
		K_\textup{mw}(x,y)
		=
		\left\{
			(z,u)\in\R\times\R^2_+\,\middle|\,
				\begin{aligned}
					&x+y\leq 1,\,-x+y\leq 1,\\
					&1-x(u_1-u_2)\leq z\leq y,\,u_1+u_2=1
				\end{aligned}
		\right\}
	\]
	which yields that, around $(\bar x,\bar y)$, the images of $K_\textup{mw}$ are uniformly bounded.
\end{example}

In \cite[Example~4.1]{LiLinZhu2024}, it has been shown that MFCQ may hold at the feasible points
of \eqref{eq:MWDref}. 
This example, however, is again awkwardly misleading as the exploited lower-level
problem does not satisfy weak Mond--Weir duality which makes the overall reformulation
\eqref{eq:MWDref} pointless. Indeed, using similar arguments as in the proof of \cref{thm:violation_of_MFCQ_Wolfe},
one can show that, in the convex situation considered here, MFCQ is violated at each feasible point of \eqref{eq:MWDref},
and in more general settings, this result remains true as long as the lower-level problem \eqref{eq:llp} satisfies
weak Mond--Weir duality for each $x\in X$.

\begin{theorem}\label{thm:violation_of_MFCQ_MondWeir}
	Let $(\bar x,\bar y,\bar z,\bar u)\in\R^n\times\R^m\times\R^m\times\R^p$ be a feasible point of \eqref{eq:MWDref}.
	Furthermore, assume that $X:=\{x\in\R^n\,|\,G(x)\leq 0\}$ holds
 	for a continuously differentiable function $G\colon\R^n\to\R^q$.
 	Then MFCQ for \eqref{eq:MWDref} is violated at $(\bar x,\bar y,\bar z,\bar u)$.
\end{theorem}

In contrast to the opinion of the authors of \cite{LiLinZhu2024},
we attest \eqref{eq:MWDref} an inherent lack of regularity at all feasible points,
and due to the presence of the surrogate variables $(z,u)$, this problem is likely to possess
more stationary points than the original bilevel optimization problem \eqref{eq:eoblp}.
Problem \eqref{eq:MWDref} shares these undesirable properties
with the reformulations \eqref{eq:KKTref}, \eqref{eq:LDref}, and \eqref{eq:WDref}
of the bilevel optimization problem \eqref{eq:eoblp}.

\subsection{Comparison in a nutshell}\label{sec:comparison}

This section presents a brief theoretical comparison of the
standard reformulations \eqref{eq:VFR}, \eqref{eq:KKTref}, and \eqref{eq:GEref},
see \cref{tab:comp_standard},
as well as the duality-based transformations \eqref{eq:LDref}, \eqref{eq:WDref}, 
and \eqref{eq:MWDref}, see \cref{tab:comp_duality_based}, 
associated with the bilevel optimization problem \eqref{eq:eoblp}.
Let us mention that a comprehensive comparison of \eqref{eq:VFR} and \eqref{eq:KKTref},
which also takes the numerical behavior of these reformulations into account, 
has already been provided in \cite{ZemkohoZhou2021}.
Most of the criteria which we exploit for our comparison here are taken from 
\cite[Table~1]{ZemkohoZhou2021}. 

\begin{table}
\centering
\begin{tabular}{l|ccc}
\toprule
 & \eqref{eq:VFR} & \eqref{eq:KKTref} & \eqref{eq:GEref} 
\\
\midrule
\# variables	&	$n+m$ & $n+m+p$ & $n+m$ 
\\
\# implicit variables & $0$ & $p$ & $0$ 
\\
\# constraints	& $p+q+1$ & $m+2p+q+1$ & $m+q$ 
\\
\midrule
lower-level convexity	& $\times$ & $\checkmark$ & $\checkmark$ 
\\
lower-level differentiability	&	$\times$	 &	$\checkmark$ & $(\checkmark)$ 
\\
lower-level regularity	&	$\times$	&	$\checkmark$ & $\times$ 
\\
\midrule
global equivalence	&	$\checkmark$ & $\checkmark$ & $\checkmark$ 
\\
local equivalence	&	$\checkmark$ & $(\times)$ & $\checkmark$ 
\\
validity of (NS-)MFCQ	&	$\times$	&	$\times$ & $-$	
\\
\bottomrule
\end{tabular}
\caption{Comparison of standand reformulations.}
\label{tab:comp_standard}
\end{table}

\begin{table}
\centering
\begin{tabular}{l|ccc}
\toprule
 &  \eqref{eq:LDref} & \eqref{eq:WDref}  	&	\eqref{eq:MWDref}
\\
\midrule
\# variables	& $n+m+p$ & $n+2m+p$ & $n+2m+p$
\\
\# implicit variables & $p$ & $ m+p$ & $m+p$
\\
\# constraints	& $2p+q+1$ & $m+2p+q+1$ & $m+2p+q+2$
\\
\midrule
lower-level convexity	& $\checkmark$ & $\checkmark$ & $\checkmark$
\\
lower-level differentiability &	$\checkmark$	&	$\checkmark$	&	$\checkmark$
\\
lower-level regularity	& $\checkmark$ &	$\checkmark$ &	$\checkmark$
\\
\midrule
global equivalence	& $\checkmark$ 	&	$\checkmark$ 	&	$\checkmark$
\\
local equivalence	& $(\times)$ & $(\times)$ &	$(\times)$
\\
validity of MFCQ & $(\times)$ & $\times$ & $\times$
\\
\bottomrule
\end{tabular}
\caption{Comparison of duality-based reformulations.}
\label{tab:comp_duality_based}
\end{table}

To start, \cref{tab:comp_standard,tab:comp_duality_based} state the number of variables
and, particularly, implicit variables (which have been made explicit)
of the different reformulations.
In fact, \eqref{eq:VFR} and \eqref{eq:GEref} do not exploit implicit variables for
the reformulation.
While \eqref{eq:KKTref} and \eqref{eq:LDref} make use of $p$ additional implicit variables,
\eqref{eq:WDref} and \eqref{eq:MWDref} require $m+p$ additional implicit variables.
Next, let us take a look at the number of scalar inequality and equality constraints 
in the different reformulations. For that purpose, we assume that 
$X=\{x\in\R^n\,|\,G(x)\leq 0\}$ holds for a function $G\colon\R^n\to\R^q$.
Furthermore, we interpret the equilibrium constraint in \eqref{eq:GEref} to be $m$-dimensional.
Let us note that 
\eqref{eq:VFR} and \eqref{eq:GEref} possess significantly less constraints than the other
four reformulations.
We have pointed out in \cref{sec:trafo_Lagrange_duality} 
that \eqref{eq:KKTref} and \eqref{eq:LDref} are closely related to each other.
However, \eqref{eq:LDref} possesses less constraints than \eqref{eq:KKTref}.

To proceed, we aim to compare requirements the lower-level problem \eqref{eq:llp}
has to satisfy in order to make the different reformulations reasonable.
A first property which we are going to investigate
is convexity of \eqref{eq:llp} w.r.t.\ the lower-level decision variable $y$
for each fixed $x\in X$. While this is essential for \eqref{eq:KKTref}, \eqref{eq:GEref},
\eqref{eq:LDref}, \eqref{eq:WDref}, and \eqref{eq:MWDref} to be reasonable,
\eqref{eq:VFR} is well-founded even in the absence of convexity.
However, as mentioned earlier, we would like to recall that convexity can be weakened 
to some generalized forms of convexity while keeping
plausibility of \eqref{eq:WDref} and \eqref{eq:MWDref}.
Reformulations \eqref{eq:KKTref}, \eqref{eq:LDref}, \eqref{eq:WDref}, and \eqref{eq:MWDref}
require continuous differentiability of the data functions $f$ and $g_1,\ldots,g_p$ 
w.r.t.\ the lower-level decision variable, 
while \eqref{eq:GEref} makes use of this property 
for the lower-level objective function $f$ only.
Let us note that whenever the lower-level problem \eqref{eq:llp} is
nonsmooth, a reformulation like \eqref{eq:KKTref} on the basis of (partial) subdifferentials 
of the involved data functions is still reasonable,
see e.g.\ \cite{DempeZemkoho2014},
and similar observations might be possible when adapting
\eqref{eq:LDref}, \eqref{eq:WDref}, and \eqref{eq:MWDref} accordingly.
It also should be noted that reformulations \eqref{eq:KKTref},
\eqref{eq:LDref}, \eqref{eq:WDref}, and \eqref{eq:MWDref}
necessarily require validity of a constraint qualification at feasible point 
of the lower-level problem, while \eqref{eq:VFR} and \eqref{eq:GEref} do not.

To finalize the comparison, 
let us have a look at the qualitative properties of the six reformulations under investigation.
First, all of them are equivalent to \eqref{eq:eoblp} in terms of global minimizers.
When local minimizers are under consideration, equivalence is preserved merely by
\eqref{eq:VFR} and \eqref{eq:GEref}.
For the other four reformulations, we only know that the local minimizers of \eqref{eq:eoblp}
can be found among the local minimizers of the respective reformulation.
The converse implication can be shown under additional assumptions 
that are not inherent but have to be checked in terms of a so-called intermediate
mapping whose properties depend on initial problem data,
see \cref{thm:LDref_local,thm:WDref_local,thm:MWDref_local} for details regarding the
duality-based reformulations.
Typically, \eqref{eq:KKTref}, \eqref{eq:LDref}, \eqref{eq:WDref}, and \eqref{eq:MWDref}
possess artificial local minimizers 
which do not correspond to local minimizers of \eqref{eq:eoblp}.
It is well known that NSMFCQ fails at all feasible points of \eqref{eq:VFR}
under mild assumptions.
Furthermore, MFCQ fails at all feasible points of \eqref{eq:KKTref} 
(as this problem possesses complementarity constraints) as well as
\eqref{eq:WDref} and \eqref{eq:MWDref}, 
see \cref{thm:violation_of_MFCQ_Wolfe,thm:violation_of_MFCQ_MondWeir}.
We have also shown in \cref{thm:violation_MFCQ_Lagrange} and the associated discussion
that NSMFCQ is likely to fail at all feasible points of \eqref{eq:LDref}.
To close, let us point out that MFCQ and NSMFCQ do not apply to the
equilibrium-constrained problem \eqref{eq:GEref} which explains the respective gap
in \cref{tab:comp_duality_based}.

	Summing up these impressions,
	from a theoretical point of view,
	the duality-based reformulations 
	\eqref{eq:LDref}, \eqref{eq:WDref}, and \eqref{eq:MWDref}
	behave similar to each other,
	and the observed behavior is close to the one of \eqref{eq:KKTref}.
	Particularly, all four reformulations suffer from the same intrinsic shortcomings.

\section{Conclusions}\label{sec:conclusions}

In this paper, we investigated transformations of the optimistic bilevel optimization problem into a single-level problem
based on the replacement of the lower-level problem by means of strong duality relations in terms of Lagrange, Wolfe,
and Mond--Weir duality. It has been shown that all three transformations suffer from the same intrinsic difficulties as 
the well-known KKT reformulation. Particularly, these single-level reformulations are likely to possess more stationary points
and, in particular, more local minimizers than the original bilevel optimization problem, and all these transformations
are quite irregular problems as a suitable version of MFCQ is, intrinsically, violated at each feasible point 
(under mild assumptions, potentially).

In the future, it has to be checked 
whether these comparatively new single-level reformulations 
can be used beneficially
in order to solve optimistic bilevel optimization problems computationally. 
Particularly, one should investigate whether these
transformations behave better than the classical KKT reformulation in numerical practice, 
as it has been indicated in
\cite{LiLinZhangZhu2022,LiLinZhu2024,OuattaraAswani2018}. 
Let us note that relaxation approaches apply to all four of
these transformations and, thus, would provide a solid base for a fair comparison, 
see e.g.\ \cite{DempeFranke2019}.


\end{document}